
\def\LaTeX{%
  \let\Begin\begin
  \let\End\end
  \let\salta\relax
  \let\finqui\relax
  \let\futuro\relax}

\def\UK{\def\our{our}\let\sz s}
\def\USA{\def\our{or}\let\sz z}

\UK



\LaTeX

\USA


\salta

\documentclass[twoside,12pt]{article}
\setlength{\textheight}{24cm}
\setlength{\textwidth}{16cm}
\setlength{\oddsidemargin}{2mm}
\setlength{\evensidemargin}{2mm}
\setlength{\topmargin}{-15mm}
\parskip2mm


\usepackage[usenames,dvipsnames]{color}
\usepackage{amsmath}
\usepackage{amsthm}
\usepackage{amssymb}
\usepackage[mathcal]{euscript}
\usepackage{cite}

%


\definecolor{viola}{rgb}{0.3,0,0.7}
\definecolor{ciclamino}{rgb}{0.5,0,0.5}

\def\gianni #1{{\color{red}#1}}
\def\juerg #1{{\color{Green}#1}}
\def\pier #1{{\color{blue}#1}}

\def\gianni #1{#1}
\def\juerg #1{#1}
\def\pier #1{#1}



\bibliographystyle{plain}


%

\finqui

\def\Beq{\Begin{equation}}
\def\Eeq{\End{equation}}
\def\Bsist{\Begin{eqnarray}}
\def\Esist{\End{eqnarray}}

\def\Bthm{\Begin{theorem}}
\def\Ethm{\End{theorem}}

\def\Bcenter{\Begin{center}}
\def\Ecenter{\End{center}}
\let\non\nonumber




\def\step #1 \par{\medskip\noindent{\bf #1.}\quad}


\def\Lip{Lip\-schitz}
\def\Holder{H\"older}

\def\aand{\quad\hbox{and}\quad}

\def\lhs{left-hand side}
\def\rhs{right-hand side}


\def\bhv{behavi\our}


\def\multibold #1{\def\arg{#1}%
  \ifx\arg\pto \let\next\relax
  \else
  \def\next{\expandafter
    \def\csname #1#1#1\endcsname{{\bf #1}}%
    \multibold}%
  \fi \next}

\def\pto{.}

\def\multical #1{\def\arg{#1}%
  \ifx\arg\pto \let\next\relax
  \else
  \def\next{\expandafter
    \def\csname cal#1\endcsname{{\cal #1}}%
    \multical}%
  \fi \next}


\def\multimathop #1 {\def\arg{#1}%
  \ifx\arg\pto\let\next\relax
  \else
  \def\next{\expandafter
    \def\csname #1\endcsname{\mathop{\rm #1}\nolimits}%
    \multimathop}%
  \fi \next}

\multibold
qwertyuiopasdfghjklzxcvbnmQWERTYUIOPASDFGHJKLZXCVBNM.

\multical
QWERTYUIOPASDFGHJKLZXCVBNM.

\multimathop
ad dist div dom meas sign supp .


\def\accorpa #1#2{\eqref{#1}--\eqref{#2}}
\def\Accorpa #1#2 #3 {\gdef #1{\eqref{#2}--\eqref{#3}}%
  \wlog{}\wlog{\string #1 -> #2 - #3}\wlog{}}


\def\separa{\noalign{\allowbreak}}

\def\neto{\mathrel{{\scriptscriptstyle\nearrow}}}
\def\seto{\mathrel{{\scriptscriptstyle\searrow}}}

\def\graffe #1{\mathopen\{#1\mathclose\}}

\def\<#1>{\mathopen\langle #1\mathclose\rangle}
\def\norma #1{\mathopen \| #1\mathclose \|}

\def\[#1]{\mathopen\langle\!\langle #1\mathclose\rangle\!\rangle}

\def\iot {\int_0^t}
\def\ioT {\int_0^T}
\def\intQt{\int_{Q_t}}
\def\intQ{\int_Q}
\def\iO{\int_\Omega}
\def\iG{\int_\Gamma}
\def\intS{\int_\Sigma}
\def\intSt{\int_{\Sigma_t}}
\def\intQi{\int_{\Qinf}}
\def\intSi{\int_{\Sinf}}

\def\dt{\partial_t}
\def\dn{\partial_\nu}

\def\cpto{\,\cdot\,}

\def\checkmmode #1{\relax\ifmmode\hbox{#1}\else{#1}\fi}
\def\aeO{\checkmmode{a.e.\ in~$\Omega$}}
\def\aeQ{\checkmmode{a.e.\ in~$Q$}}

\def\AeQ{\checkmmode{a.e.\ in~$\Qinf$}}

\def\Aet{\checkmmode{a.e.\ in~$(0,+\infty)$}}

\def\limn{\lim_{n\to\infty}}


\def\erre{{\mathbb{R}}}




\def\genspazio #1#2#3#4#5{#1^{#2}(#5,#4;#3)}
\def\spazio #1#2#3{\genspazio {#1}{#2}{#3}T0}

\def\L {\spazio L}
\def\H {\spazio H}
\def\W {\spazio W}

\def\C #1#2{C^{#1}([0,T];#2)}
\def\spazioinf #1#2#3{\genspazio {#1}{#2}{#3}{+\infty}0}
\def\LL {\spazioinf L}


\def\Lx #1{L^{#1}(\Omega)}
\def\Hx #1{H^{#1}(\Omega)}
\def\Wx #1{W^{#1}(\Omega)}
\def\LxG #1{L^{#1}(\Gamma)}
\def\HxG #1{H^{#1}(\Gamma)}

\def\Cx #1{C^{#1}(\overline\Omega)}
\def\CxG #1{C^{#1}(\Gamma)}
\def\LQ #1{L^{#1}(Q)}

\def\Luno{\Lx 1}
\def\Ldue{\Lx 2}

\def\Huno{\Hx 1}
\def\Hdue{\Hx 2}

\def\HunoG{\HxG 1}
\def\HdueG{\HxG 2}

\def\LunoG{\LxG 1}
\def\LdueG{\LxG 2}


\def\LQ #1{L^{#1}(Q)}


\let\theta\vartheta
\let\eps\varepsilon
\let\phi\varphi

\let\TeXchi\chi                         
\newbox\chibox
\setbox0 \hbox{\mathsurround0pt $\TeXchi$}
\setbox\chibox \hbox{\raise\dp0 \box 0 }
\def\chi{\copy\chibox}


\def\QED{\hfill $\square$}


\def\longtime{long-time}
\def\Longtime{Long-time}
\def\omegalimit{$\omega$-limit}
\def\CO{C_\Omega}

\def\Qinf{Q_\infty}
\def\Sinf{\Sigma_\infty}

\def\suG{_{|\Gamma}}

\def\VG{V_\Gamma}
\def\HG{H_\Gamma}

\def\nablaG{\nabla_\Gamma}
\def\DeltaG{\Delta_\Gamma}
\def\rhoG{\rho_\Gamma}
\def\FG{F_\Gamma}
\def\vG{v_\Gamma}
\def\zetaG{\zeta_\Gamma}

\def\muz{\mu_0}
\def\rhoz{\rho_0}
\def\rhozG{{\rhoz}\suG}

\def\mueps{\mu^\eps}
\def\rhoeps{\rho^\eps}
\def\rhoGeps{\rho_\Gamma^\eps}

\def\mus{\overline\mu}
\def\rhos{\overline\rho}
\def\rhoGs{\overline\rho_\Gamma}
\def\mus{\mu_s}
\def\rhos{\rho_s}
\def\rhoGs{\rho_\juerg{{s_\Gamma}}}

\def\muo{\mu_\omega}
\def\rhoo{\rho_\omega}
\def\rhoGo{\rho_\juerg{{\omega_\Gamma}}}

\def\mui{\mu^\infty}
\def\rhoi{\rho^\infty}
\def\rhoGi{\rhoG^\infty}

\def\tn{t_n}
\def\mun{\mu^n}
\def\rhon{\rho^n}
\def\rhoGn{\rho_\Gamma^n}

\def\rhomin{\rho_*}
\def\rhomax{\rho^*}
\def\gstar{g_*}

\def\Vp{V^*}
\def\normaV #1{\norma{#1}_V}
\def\normaH #1{\norma{#1}_H}
\def\normaW #1{\norma{#1}_W}
\def\normaVG #1{\norma{#1}_{\VG}}

\let\hat\widehat

\Begin{document}


%
\title{Limiting problems\\[0.3cm] 
  for a nonstandard viscous Cahn--Hilliard system\\[0.3cm] 
  with dynamic boundary conditions}
\author{}
\date{}
\maketitle
\Bcenter
\vskip-1cm
{\large\sc Pierluigi Colli$^{(1)}$}\\
{\normalsize e-mail: {\tt pierluigi.colli@unipv.it}}\\[.25cm]
{\large\sc Gianni Gilardi$^{(1)}$}\\
{\normalsize e-mail: {\tt gianni.gilardi@unipv.it}}\\[.25cm]
{\large\sc J\"urgen Sprekels$^{(2)}$}\\
{\normalsize e-mail: {\tt sprekels@wias-berlin.de}}\\[.45cm]
$^{(1)}$
{\small Dipartimento di Matematica ``F. Casorati'', Universit\`a di Pavia}\\
{\small via Ferrata 5, 27100 Pavia, Italy}\\[.3cm]
$^{(2)}$
{\small Department of Mathematics}\\
{\small Humboldt-Universit\"at zu Berlin}\\
{\small Unter den Linden 6, 10099 Berlin, Germany}\\[2mm]
{\small and}\\[2mm]
{\small Weierstrass Institute for Applied Analysis and Stochastics}\\
{\small Mohrenstrasse 39, 10117 Berlin, Germany}
\Ecenter
\Begin{abstract}\noindent
This note is concerned with a nonlinear diffusion problem of phase-field type,
consisting of a parabolic system of two partial differential equations, complemented by boundary and initial conditions. The system arises from a
model of two-species phase segregation on an atomic lattice and was introduced 
by Podio-Guidugli in Ric.\ Mat.\ {\bf 55} (2006), pp.~105--118. The two unknowns are the phase parameter and the chemical potential. In contrast to previous investigations about this PDE system, we consider here a dynamic boundary condition for the phase variable that involves the Laplace-Beltrami operator and models an additional nonconserving phase transition occurring on the surface of the domain. We are interested to some asymptotic analysis and first discuss the asymptotic limit of the system as the viscosity coefficient of the order parameter equation tends to 0: the  convergence of solutions to the corresponding solutions for the limit problem is proven. Then, we study the long-time behavior of the system for both problems, with positive or zero viscosity coefficient, and characterize the omega-limit set in both cases. 
\\[4mm]
\End{abstract}

\salta
\pagestyle{myheadings}
\newcommand\testopari{\sc Colli \ --- \ Gilardi \ --- \ Sprekels}
\newcommand\testodispari{\sc \juerg{Limits for a 
nonstandard Cahn--Hilliard system with dynamic b.c.}}
\markboth{\testodispari}{\testopari}
\finqui
%

\section{Introduction}
\label{Intro}
\setcounter{equation}{0}

A recent line of research originated from the following evolutionary system
of partial differential equations:
\Bsist
  && 2\rho \, \dt\mu
  + \mu \, \dt\rho
  - \Delta\mu = 0
  \aand
  \mu \geq 0
  \label{oldprima}
  \\
  && - \Delta\rho + F'(\rho) = \mu
   \label{oldseconda}
\Esist
in $\Qinf:=\Omega\times(0,+\infty)$,
where $\Omega \subset \erre^3$ is a bounded and smooth domain with boundary~$\Gamma$. 
The system \eqref{oldprima}--\eqref{oldseconda} comes out from a model for phase segregation through atom rearrangement on a lattice that has been proposed by Podio-Guidugli~\cite{PG}. This model (see also \cite{CGPS3} for a detailed derivation) is a modification of the Fried--Gurtin approach to
phase segregation processes (cf.~\cite{FG,G}). The order parameter 
$\rho$, which in many cases represents the (normalized) density of one of the phases, 
and the chemical potential $\mu$ are the unknowns of the system.
Moreover, $F'$ represents the derivative of a double-well potential~$F$.
Besides everywhere defined potentials, a~typical and important example of $F$ 
is the so--called {\em logarithmic double-well potential\/} given~by
\Beq
  F_{log}(r) := (1+r)\ln (1+r)+(1-r)\ln (1-r) + \alpha_1  (1- r^2) + \alpha_2 r,
  \quad r \in (-1,1),
  \label{logpot}
\Eeq
for some real coefficients $\alpha_1, \,  \alpha_2 $. Note that, if $\alpha_2$ is taken null and $\alpha_1 > 1$, it turns out that $F$ actually exhibits two wells, with a local maximum at $r=0$. In the case when $\alpha_2 \not=0$, then one of the two minima of $F$ is preferred, in the sense that there is a global minimum point (positive if $\alpha_2 <0$, negative if $\alpha_2 >0$) of the function. As a particular feature of 
\eqref{logpot}, observe that the derivative of the logarithmic potential becomes singular at $\,\pm 1$.

About equations \eqref{oldprima} and \eqref{oldseconda}, we point out that
the model developed in~\cite{PG} is based on a local 
free energy density (in the bulk) of the form
\begin{equation}
\label{fe1}
\psi(\rho,\nabla\rho,\mu)=-\mu\,\rho+F(\rho)+\frac 1 2\,|\nabla\rho|^2,
\end{equation}
so to derive equations \eqref{oldprima}--\eqref{oldseconda}, which must be 
complemented with boundary and initial conditions.
As \juerg{far as} the former are concerned, the standard boundary conditions for this class of problems are the homogeneous Neumann ones, namely
\Beq
  \dn\mu = \dn\rho = 0 
  \quad \hbox{on \,$\Sinf:=\Gamma\times(0,+\infty)$},
  \label{IneumannBC}
\Eeq
where $\dn$ denotes the outward normal derivative. Combining now 
\eqref{oldprima}--\eqref{oldseconda} with \eqref{IneumannBC}, we obtain a 
set of equations and conditions that is a variation of the celebrated 
Cahn--Hilliard system originally introduced in \cite{CahH}
and first studied mathematically in \cite{EllSh} (for an updated
list of references on the Cahn--Hilliard system, see~\cite{Heida}). Nonetheless, 
an initial value problem for \eqref{oldprima}--\eqref{oldseconda}, \eqref{IneumannBC} turns out to be strongly ill-posed (see \cite[Subsect.~1.4]{CGPS8}, where an example is given): indeed, the related problem may have infinitely many smooth and even nonsmooth solutions. 
Then, two small regularizing parameters $\varepsilon>0$ and $\delta>0$ were introduced and considered in \cite{CGPS3}, which led to the regularized model equations 
\Bsist
  && \bigl( \eps + 2\rho \bigr) \, \dt\mu
  + \mu \, \dt\rho
  - \Delta\mu = 0 \,,
  \label{oldprimaeps}
  \\
  && \delta\, \dt\rho - \Delta\rho + F'(\rho) = \mu\,. 
   \label{oldsecondaeps}
\Esist
This regularized system has been deeply examined in \cite{CGPS3}, when both $\eps$ and 
$\delta$ are positive and fixed. In addition, let us underline that, 
while one can let $\eps$ tend to zero (see~\cite{CGPS4}) 
and obtain a solution to the limiting problem with $\eps=0$,
it seems extremely difficult to the pass to the limit as $\delta$ goes to $0$. 
In fact, ill-posedness still holds for $\delta=0$, even if $\eps$ is kept positive.
Hence, one has to assume that $\delta$ is a fixed positive coefficient.
Therefore, from now on, we take $\delta=1$, without loss of generality.
Let us point out that the \longtime\ \bhv\ of the solutions has been studied both with $\eps>0$ (cf.~\cite{CGPS3}) and $\eps=0$ (cf.~\cite{CGPS4}).

The system \eqref{oldprimaeps}--\eqref{oldsecondaeps} constitutes
a modification of the so-called {\em viscous} Cahn--Hill\-iard system (see \cite{NC} and the recent contributions\cite{CGW, CGS0, CGS2} along with their references). 
We point out that \eqref{oldprimaeps}--\eqref{oldsecondaeps} 
was analyzed, in the case of the
boundary conditions \eqref{IneumannBC}, in the 
papers \cite{CGPS3,  CGPSco, CGSco1} concerning
well-posedness, regularity, and optimal control. Later, the
local free energy density \eqref{fe1} was generalized to the form
\begin{equation}
\label{fe2}
\psi(\rho,\nabla\rho,\mu)=-\mu\,g(\rho)+F(\rho)+\frac 1 2\,|\nabla\rho|^2,
\end{equation}
thus putting $g(\rho)$ in place of~$\rho$,
where $g$ is a nonnegative function on the domain of~$F$.
This leads to the system
\Bsist
  && \bigl( \eps + 2g(\rho) \bigr) \, \dt\mu
  + \mu \, g'(\rho) \, \dt\rho
  - \Delta\mu = 0,
  \label{Igenprimaeps}
  \\
  && \dt\rho - \Delta\rho + F'(\rho) = \mu \,g'(\rho),
   \label{Igensecondaeps}
\Esist
which is a generalization of \accorpa{oldprimaeps}{oldsecondaeps}
and has been studied in \cite{CGPS7,CGPS6} \juerg{for the case} $\eps=1$. Let us mention also the contribution 
\cite{CGKPS} dealing with the time discretization of the problem and proving convergence results and error estimates. The related phase relaxation system
(in which the diffusive term $- \Delta\rho $ disappears from \eqref{Igensecondaeps}), has been dealt with in \cite{CGKS1, CGKS2, CGS-lom}. 
We also point out the recent papers \cite{CGS3, CGS4, CGS5},
where a nonlocal version of \eqref{Igenprimaeps}--\eqref{Igensecondaeps} -- based on the replacement of the diffusive term of \eqref{Igensecondaeps} with a nonlocal operator acting on $\rho$ -- has been largely investigated, also from the side of optimal control.

Now, if we take $\eps=0$ in \eqref{Igenprimaeps}--\eqref{Igensecondaeps}, we obtain
\Bsist
  && 2g(\rho) \, \dt\mu
  + \mu \, g'(\rho) \, \dt\rho
  - \Delta\mu = 0
  \label{Iprima}
  \\
  && \dt\rho - \Delta\rho + F'(\rho) = \mu \,g'(\rho),
   \label{Iseconda}
\Esist
which looks \juerg{like a generalization} of the viscous version of~\accorpa{oldprima}{oldseconda}, where the affine function $\rho \mapsto \rho$ is replaced by a concave function $\rho \mapsto g(\rho)$, with $g$ possessing suitable properties that are made precise in the later assumption~\eqref{hpg}. In particular, the new $g$ may be symmetric and strictly concave: a possible simple choice of $g$ satisfying \eqref{hpg} is 
\Beq
  g(r) = 1- r^2 , 
  \quad r \in [-1,1]. 
  \label{gpot}
\Eeq
Note that, if one collects \eqref{logpot} and \eqref{gpot} and assumes $\alpha_2\not=0$, the combined function
\begin{equation}
\label{combi}
-  \mu g(\rho) + F_{log} (\rho)   \quad \hbox{(which is a part of }\, \psi) 
\end{equation}
shows a global minimum in all cases, and it depends on the values of $(\alpha_1 - \mu)$ and $\alpha_2$ which minimum actually occurs. Let us notice that the function 
in \eqref{combi} turns out to be convex in the whole of $(-1,1)$ for sufficiently large values of $\mu$. On the other hand, the framework fixed by assumptions \eqref{hpg}--\eqref{domination} allows for more general choices of $g$ and $F $.    

However, until now we remain at the situation that the boundary conditions are of Neumann type \juerg{for both} $\mu$ and~$\rho$, as in~\eqref{IneumannBC}. Instead,
in the present work we  treat the dynamic boundary condition for~$\rho$,
i.e., we complement the above systems with
\Beq
  \dn\mu = 0
  \aand
  \dn\rho + \dt\rhoG - \DeltaG\rhoG + \FG'(\rhoG) = 0
  \quad \hbox{on $\Sinf$},
  \label{IdynBC}
\Eeq
where $\rhoG$ is the trace of~$\rho$,
$\DeltaG$~is the Laplace-Beltrami operator on the boundary,
$\FG'$~is the derivative of another potential~$\FG$ having more or less the same 
behavior as $F$,
and the \rhs\ of the dynamic boundary condition \juerg{equals} zero, just for simplicity.
Indeed, one could consider a \juerg{nonzero} forcing term
satisfying proper assumptions, as done in~\cite{CGS10}. 
Once again, we have to add initial conditions. 

Thus, we are concerned with a total free energy of the system which also includes 
a contribution on the boundary; in fact, we postulate that a phase transition phenomenon is occurring as well on the boundary, and the physical variable on the boundary is just the trace of the phase variable in the bulk. This corresponds to
a total free energy functional of the form
\begin{align}
\label{fe3}
\mathbf{\Psi} [\rho\juerg{(t)}, \rho_\Gamma\juerg{(t)},\mu(t)]
= &\int_\Omega \Bigl[-\mu\,g(\rho)+F(\rho)+\frac 1 2\,|\nabla\rho|^2\Bigr]
\juerg{(t)} \non\\
&+  \int_\Gamma \Bigl[[-u_\Gamma\,\rho_\Gamma + F_\Gamma (\rho_\Gamma)+\frac 1 2\,|\nabla_\Gamma \rho_\Gamma|^2\Bigr]\juerg{(t)}, \quad t\geq0,
\end{align}
where $\nabla_\Gamma$ is the surface gradient and $u_\Gamma$ may stand for the source term that exerts a (boundary) control on the system. From this expression of the total free energy, one recovers the PDE system resulting from equations 
\eqref{Iprima}--\eqref{Iseconda} and the boundary conditions \eqref{IdynBC},
with $u_\Gamma$ in place of $0$ in the \rhs\ of the second condition.
In relation to this, we would like to mention the contribution \cite{CGS11} 
dealing with the optimal boundary control problem for 
the system \accorpa{oldprimaeps}{oldsecondaeps}, \eqref{IdynBC} with $\eps=1$.

As for the dynamic boundary conditions, we would like to add some comments 
on the recent growing interest in the mathematical literature,
either for the justification (see, e.g., \cite{FiMD1, FiMD2, Li}) 
or for the investigation of systems including dynamic boundary conditions. 
Without trying to be exhaustive, we point out at least the 
contributions \cite{Calcol, ChGaMi, ChFaPr, CF1, CF2, CF3, CGS0, CGS1, CGS2, CS, CGM1, CGM2, GaGu, GaWa, GiMiSc, GMS10, GoMi, GMS11, Is, MRSS, MiZe, PrRaZh, RaZh},
which are concerned with various types of systems endowed with the dynamic boundary conditions for either some or all of the unknowns. Our citations mostly refer to 
phase-field models involving the Allen--Cahn and Cahn--Hilliard equations, 
whose structure is generally simpler than the one considered in the present paper.

Our aim here is investigating the \longtime\ \bhv\ of the full system
in both \juerg{the} cases $\eps>0$ and $\eps=0$
(similarly as in \cite{CGPS3,CGPS4}, in which the Neumann 
boundary conditions \eqref{IneumannBC} were considered).
More precisely, we show that the \omegalimit\ of any trajectory in a suitable topology
\juerg{consists only} of stationary solutions.
In order to treat this problem also with $\eps=0$,
we first study the asymptotics as $\eps$ tends to zero.
To do that, we underline that the reasonable and somehow natural assumptions \eqref{hpg} for $g$ along with the requirements \eqref{hpFprimo}--\eqref{domination} on  $F$ and $\FG$ allow us to show that the variables $\rho$ and $\rho_\Gamma$ are strictly separated from the (singular) values $\pm 1$. Indeed, we can prove this separation property and obtain the strict positivity of $g(\rho)$ as a consequence; moreover, we believe that our approach here is physically reasonable and more satisfactory than the analysis and results obtained in~\cite{CGPS4}. 

The paper is organized as follows: 
in the next section, we list our assumptions and notations and state our results,
while the corresponding proofs are given in the last two sections.
Precisely, in Section~\ref{WELLPOSEDNESS}, 
we perform the asymptotic analysis as $\eps$ tends to zero
and prove the well-posedness of the problem for $\eps=0$;
in Section~\ref{LONGTIME}, we study the \longtime\ \bhv\ of the solution 
\juerg{under the assumption} $\eps\geq0$.


\section{Statement of the problem and results}
\label{STATEMENT}
\setcounter{equation}{0}

In this section, we state precise assumptions and notations and present our results.
First of all, the \juerg{set} $\Omega\subset\erre^3$ is assumed to be bounded, connected and smooth.
As in the Introduction, $\dn$~and $\DeltaG$ stand for the \juerg{outward} normal derivative
and the Laplace-Beltrami operator on the boundary~$\Gamma$.
Furthermore, we denote by $\nablaG$ the surface gradient.

If $X$ is a (real) Banach space, $\norma\cpto_X$ denotes both its norm and the norm of~$X^3$,
$X^*$~is its dual space, and ${}_{X^*}\<\cpto,\cpto>_X$ is the dual pairing between $X^*$ and~$X$. 
The only exception from this convention is given
by the $L^p$ spaces, $1\leq p\leq\infty$, for which we use the abbreviating notation
$\norma\cpto_p$ for the norms in~$L^p(\Omega)$. 
Furthermore, we put
\Bsist
  && H := \Ldue \,, \quad  
  V := \Huno 
  \aand
  W := \graffe{v\in\Hdue : \ \dn v=0},
  \label{defspaziO}
  \\
  && \HG := \LdueG 
    \aand
  \VG := \HunoG ,
  \label{defspaziG}
  \\
  && \calH := H \times \HG 
  \aand
  \calV := \graffe{(v,\vG) \in V \times \VG : \ \vG = v\suG}.
  \label{defspaziprod}
\Esist
\Accorpa\Defspazi defspaziO defspaziprod
We also set, for convenience,
\gianni{%
\Bsist
  && Q_t := \Omega \times (0,t)
  \aand
  \Sigma_t := \Gamma \times (0,t)
  \quad \hbox{for $0<t<+\infty$},
  \non
  \\
  && \Qinf := \Omega\times (0,+\infty)
  \aand
  \Sinf := \Gamma\times (0,+\infty),
  \label{defQt}
\Esist
}%
and often use the shorter notations $Q$ and $\Sigma$ \,if\, $t=T$, a fixed final time~$T\in(0,+\infty)$.

Now, we list our assumptions.
For the structure of our system, we are given three functions
$g\in C^2[-1,1]$ and $F,\,\FG\in C^2(-1,1)$ which satisfy
\Bsist
  && g \geq 0 , \quad
  g'' \leq 0 , \quad
  g'(-1) > 0 
  \aand
  g'(1) < 0,
  \label{hpg}
  \\
  \separa
  && \lim_{r\seto -1}F'(r)
  = \lim_{r\seto -1}\FG'(r)
  = -\infty
  \aand
  \lim_{r\neto 1}F'(r)
  = \lim_{r\neto 1}\FG'(r)
  = +\infty, 
  \qquad
  \label{hpFprimo}
  \\
  \separa
  && F''(r) \geq - C
  \aand
  \FG''(r) \geq -C,
  \quad \hbox{for every $r\in(-1,1)$},
  \label{hpFsecondo}
  \\
  \separa
  && |F'(r)| \leq \eta |\FG'(r)| + C
  \quad \hbox{for every $r\in(-1,1)$},
  \label{domination}
\Esist
\Accorpa\HPstruttura hpg domination
\juerg{with} some positive constants $C$ and~$\eta$.

For the initial data, we make rather strong assumptions
in order to apply the results of~\cite{CGS10} without any trouble.
However, \juerg{our first assumption on $\muz$ could} be replaced by $\muz\in V$.
Precisely, we assume that
\begin{align}
  & \muz \in W 
  \aand
  \muz \geq 0 \quad \hbox{in $\Omega$}\,\pier{;}
  \label{hpmuz}
  \\
  & \rhoz  \in \Hdue \,, \quad
  \rhozG \in \HdueG \,, \quad
  \min\rhoz > -1
  \aand
  \max\rhoz < 1 \,.
  \label{hprhoz}
\end{align}
\Accorpa\HPdati hpmuz hprhoz

At this point, we are ready to state our problem.
For $\eps\geq0$, we look for a triplet $(\mu,\rho,\rhoG)$ satisfying  
the regularity requirements and solving the problem stated below.
As for the regularity, we pretend that
\begin{align}
  & \mu \in \H1H \cap \C0V \cap \L2W, 
  \label{regmu}
  \\
  & (\rho,\rhoG) \in \W{1,\infty}\calH \cap \H1\calV \cap \L\infty{\Hdue\times\HdueG},
  \label{regrho}
  \\
  & \mu \geq 0 \,, \quad
  -1 < \rho < 1
  \aand
  (F'(\rho),\FG'(\rhoG)) \in \L\infty\calH , 
  \label{disugsoluz}
\end{align}
\Accorpa\Regsoluz regmu disugsoluz
for every finite $T>0$, and the problem \juerg{reads}
\begin{align}
  & \bigl( \eps + 2g(\rho) \bigr) \dt\mu
  + \mu g'(\rho) \dt\rho
  - \Delta\mu
  = 0
  \quad \AeQ\,,
  \label{prima}
  \\
  & \iO \dt\rho \, v
  + \iG \dt\rhoG \, \vG
  + \iO \nabla\rho \cdot \nabla v
  + \iG \nablaG\rhoG \cdot \nablaG\vG
  + \iO F'(\rho) v
  + \iG \FG'(\rhoG) \vG 
  \non
  \\
  & = \iO \mu g'(\rho) v
  \qquad \hbox{\Aet\ and for every $(v,\vG)\in\calV$}\,,
  \label{seconda}
  \\
  & \mu(0) = \muz
  \aand
  \rho(0) = \rhoz
  \quad \aeO \,.
  \label{cauchy}
\end{align}
\Accorpa\Pbl prima cauchy
Notice that the Neumann boundary condition $\dn\mu=0$ and
the fact that $\rhoG$ is the trace of $\rho$ on~$\Sigma$ 
are contained in \eqref{regmu} and~\eqref{regrho}, respectively,
due to the definitions \Defspazi\ of the spaces involved.
By accounting \juerg{for} the regularity conditions~\Regsoluz,
it is clear that the variational problem~\eqref{seconda} is equivalent~to
\begin{align} 
  & \dt\rho - \Delta\rho + F'(\rho) = \mu g'(\rho)
  \quad \hbox{in $\Qinf$}\,,
  \label{secondaO}
  \\
  & \dn\rho + \dt\rhoG - \DeltaG\rhoG + \FG'(\rhoG) = 0
  \quad \hbox{on $\Sinf$}\,.
  \label{secondaG}
\end{align}
\Accorpa\Secondabvp secondaO secondaG
\juerg{Moreover, it follows from standard embedding results (see, e.g., \cite[Sect.\,8, Cor.\,4]{Simon})
that $\rho\in C^0(\overline Q)$ and thus also $\rhoG\in C^0(\overline\Sigma)$.}

Our starting point is the well-posedness result for $\eps>0$ that we state below
and is already known.
Indeed, recalling \accorpa{hpFprimo}{hpFsecondo}, we set
\Beq
  \hat\beta(r) := F(r) - F(0) - F'(0) r + \frac C2 \, r^2
  \quad \hbox{for $r\in(-1,1)$}
  \aand
  \hat\pi := F- \hat\beta,
  \non
\Eeq
and analogously introduce $\hat\beta_\Gamma$ and $\hat\pi_\Gamma$, starting from~$\FG$.
Then, we consider the convex and lower semicontinuous extensions
of  $\hat\beta$ and $\hat\beta_\Gamma$ to the whole of $\erre$
and smooth extensions of $\hat\pi$ and $\hat\pi_\Gamma$
with bounded second derivatives.
Therefore, the assumptions of \cite[Thm.~2.1]{CGS10} are satisfied
and the following \juerg{well-posedness} result holds true.

\Bthm
\label{Wellposednesseps}
Assume \HPstruttura\ and $\eps>0$ for the structure and \HPdati\ for the initial data.
Then problem \Pbl\ has a unique solution $(\mueps,\rhoeps,\rhoGeps)$ satisfying the regularity properties \Regsoluz.
\Ethm

Our aim is the following:
$i)$~by starting from the solution $(\mueps,\rhoeps,\rhoGeps)$,
we let $\eps$ tend to zero and prove that problem \Pbl\ with $\eps=0$ has a solution $(\mu,\rho,\rhoG)$;
$ii)$~such a solution is unique;
$iii)$~for $\eps\geq0$, we study the \omegalimit\ of every trajectory.

Indeed, for $i)$ and~$ii)$, we prove the following result in Section~\ref{WELLPOSEDNESS}:

\Bthm
\label{Wellposedness}
Assume \HPstruttura\ for the structure and \HPdati\ for the initial data.
Then problem \Pbl\ with $\eps=0$ has a unique solution $(\mu,\rho,\rhoG)$ 
satisfying the regularity properties \Regsoluz.
Moreover, for some constants $\rhomin,\rhomax\in(-1,1)$
that depend only on the shape of the nonlinearities and on the initial data,
both $(\mu,\rho,\rhoG)$ 
and the solution $(\mueps,\rhoeps,\rhoGeps)$ given by Theorem~\ref{Wellposednesseps} 
satisfy the separation property
\Beq
  \rhomin \leq \rho \leq \rhomax
  \aand
  \rhomin \leq \rhoeps \leq \rhomax
  \quad \hbox{in \juerg{$\,\overline\Omega\times [0,+\infty)$}}.
  \label{separaz}
\Eeq
Finally, $(\mueps,\rhoeps,\rhoGeps)$
converges to $(\mu,\rho,\rhoG)$ in a proper topology.
\Ethm

The last Section~\ref{LONGTIME} is devoted to study the \longtime\ \bhv\ of the solution
in both \juerg{the} cases $\eps>0$ and $\eps=0$.
To this end, for a fixed $\eps\geq0$, 
we use the simpler symbol $(\mu,\rho,\rhoG)$ for the solution on $[0,+\infty)$
and observe that the regularity~\Regsoluz\ on every finite time interval
implies that $(\mu,\rho,\rhoG)$ is a continuous $(H\times\calV)$-valued function.
In particular, it can be evaluated at every time~$t$,
and the following definition of \omegalimit\ is completely meaningful:
\Bsist
  && \omega(\mu,\rho,\rhoG)
  := \Bigl\{
    (\muo,\rhoo,\rhoGo) \in H \times \calV :\ 
    (\mu,\rho,\rhoG)(\tn) \to (\muo,\rhoo,\rhoGo)
    \qquad
  \non
  \\
  && \phantom{\omega(\mu,\rho,\rhoG)
  := \Bigl\{}
  \hbox{weakly in $H\times\calV$ for some sequence $\tn\nearrow+\infty$}
  \Bigr\}
  \label{omegalimit}.
\Esist
Besides, we consider the stationary solutions.
It is immediately seen that a stationary solution is a
triplet $(\mus,\rhos,\rhoGs)$ satisfying the following conditions:
the first component $\mus$ is a constant,
and $(\rhos,\rhoGs)\in\calV$ is a solution to the system
\Bsist
  && \iO \nabla\rhos \cdot \nabla v
  + \iG \nablaG\rhoGs \cdot \nablaG\vG
  + \iO F'(\rhos) v
  + \iG \juerg{F'_\Gamma(\rhoGs)} \vG
  \non
  \\
  && = \iO \mus \,g'(\rhos) v
  \qquad \hbox{for every $(v,\vG)\in\calV$}.
  \label{secondastaz}
\Esist
In terms of a boundary value problem, the conditions $(\rhos,\rhoGs)\in\calV$ and \eqref{secondastaz} mean
that 
\Beq  
  - \Delta\rhos + F'(\rhos) = \mus \,g'(\rhos)
  \quad \hbox{in $\Omega$},
  \quad
  \rhoGs = {\rhos}\suG
  \aand
  \dn\rhos - \DeltaG\rhoGs + \FG'(\rhoGs) = 0
  \quad \hbox{on $\Gamma$}.
  \non
\Eeq
We prove the following result:

\Bthm
\label{Omegalimit}
Assume \HPstruttura\ and $\eps\geq0$ for the structure and \HPdati\ for the initial data,
and let $(\mu,\rho,\rhoG)$ \juerg{be} the unique solution to problem \Pbl\
satisfying the regularity requirements \Regsoluz.
Then the \omegalimit\ \eqref{omegalimit} is nonempty and consists only of stationary solutions.
In particular, there exists a constant $\mus$ such that
problem \eqref{secondastaz} has at least \juerg{one} solution $(\rhos,\rhoGs)\in\calV$.
\Ethm

Throughout the paper, we will repeatedly use the Young inequality
\Beq
  a\,b \leq \delta\,a^2 + \frac 1{4\delta} \, b^2
  \quad \hbox{for all $a,b\in\erre$ and $\delta>0$},
  \label{young}
\Eeq
as well as the \Holder\ inequality and the continuity of
the embedding $V\subset L^p(\Omega)$ for every $p\in[1,6]$
(since $\Omega$ is \juerg{three-dimensional}, bounded and smooth).
\juerg{Besides, this embedding is compact for $p<6$, 
and also the embedding $W\subset\Cx0$ is compact. In particular, we have the 
compactness inequality
\begin{equation}
\label{compact}
\|v\|_4\,\le\,\delta\,\|\nabla v\|_2\,+\,\widetilde C_\delta\,\|v\|_2
\quad\mbox{for every $v\in H^1(\Omega)$ and $\delta>0$},
\end{equation}
where $\widetilde C_\delta$ depends only on $\Omega$ and $\delta$.}
We also recall some well-known estimates \juerg{from} trace theory
and \juerg{from} the theory of elliptic equations we use in the sequel.
For \juerg{any} $v$ and $\vG$ that make the \rhs s meaningful,
we have \juerg{that}
\begin{align}
  & \norma{\dn v}_{\HxG{-1/2}}
  \leq \CO \bigl( \norma v_{\Huno} + \norma{\Delta v}_{\Ldue} \bigr)\,, 
  \label{dnHmum}
  \\[1mm]
  & \norma{\dn v}_{\LdueG}
  \leq \CO \bigl( \norma v_{\Hx{3/2}} + \norma{\Delta v}_{\Ldue} \bigr)\,,
  \label{dnLd}
  \\[1mm]
  & \norma v_{\Hdue} \leq \CO \bigl( \norma{v\suG}_{\HxG{3/2}} + \norma{\Delta v}_{\Ldue} \bigr)\,,
  \label{stimaHdueD}
  \\[1mm]
  & \norma v_{\Hdue} \leq \CO \bigl( \norma v_{\Huno} + \norma{\Delta v}_{\Ldue} \bigr)
  \,\quad \hbox{if $\dn v=0$ \ \pier{on $\Gamma$}}\,,
  \label{stimaHdueN}
  \\[1mm]
  & \norma\vG_{\HdueG} \leq \CO \bigl( \norma\vG_{\HunoG} + \norma{\DeltaG\vG}_{\LdueG} \bigr)\,,
  \label{stimaHdueG}
  \\[1mm]
  & \norma\vG_{\HxG{3/2}} \leq \CO \bigl( \norma\vG_{\HunoG} + \norma{\DeltaG\vG}_{\HxG{-1/2}} \bigr)\,,
  \label{stimaHtmG}
\end{align}
\Accorpa\Richiami dnHmum stimaHtmG
\juerg{with a constant $\CO>0$ that} depends only on~$\Omega$.

We conclude this section by stating a general rule concerning the constants 
that appear in the estimates \juerg{to be performed} in the sequel.
The small-case symbol $c$ stands for a generic constant
whose values might change from line to line and even \juerg{within} the same line
and \juerg{depends} only on~$\Omega$, on the shape of the nonlinearities,
and on the constants and the norms of the functions involved in the assumptions of our statements.
In particular, the values of $c$ do not depend on $\eps$ and~$T$ if the latter is considered.
A~small-case symbol with a subscript like $c_\delta$ (in~particular, with $\delta=T$)
indicates that the constant might depend on the parameter~$\delta$, in addition.
On the contrary, we mark precise constants that we can refer~to
by using different symbols, like in~\accorpa{hpFsecondo}{domination} and 
\eqref{compact}--\eqref{stimaHtmG}.


\section{Well-posedness}
\label{WELLPOSEDNESS}
\setcounter{equation}{0}

This section is devoted to the proof of Theorem~\ref{Wellposedness}.
First, we prove the separation properties~\eqref{separaz}.
Then, we show uniqueness.
Finally, we prove convergence for the family $\graffe{(\mueps,\rhoeps,\rhoGeps)}$
and derive existence for the problem with $\eps=0$.

\step
Separation

We assume that $\eps\geq0$ and that $(\mu,\rho,\rhoG)$
is a solution to problem \Pbl\ satisfying \Regsoluz.
Recalling \eqref{hprhoz} and \accorpa{hpg}{hpFsecondo},
we \juerg{may} choose $\rhomin,\rhomax\in(-1,1)$ such that
$\rhomin\leq\rhoz\leq\rhomax$ \juerg{and}
\Bsist
  &&g'(r) > 0 
  \aand
  F'(r) < 0
  \quad \hbox{for $-1<r\leq\rhomin$},
  \non
  \\
  &&g'(r) < 0 
  \aand
  F'(r) > 0
  \quad \hbox{for $\rhomax\leq r<1$}.
  \non
\Esist
Now, we show that $\rhomin\leq\rho\leq\rhomax$, using the positivity of~$\mu$ (see~\eqref{disugsoluz}).
In fact, we prove just the upper inequality, since the \juerg{proof of the other} is similar.
We test \eqref{seconda}, written at the time~$s$, by $((\rho-\rhomax)^+,(\rhoG-\rhomax)^+)(s)$
and integrate over~$(0,t)$ with respect to~$s$.
We have
\Bsist
  && \frac 12 \iO |(\rho(t)-\rhomax)\juerg{^+}|^2
  + \frac 12 \iG |(\rhoG(t)-\rhomax)\juerg{^+}|^2
  + \intQt |\nabla(\rho-\rhomax)\juerg{^+}|^2
  + \intSt |\nablaG(\rhoG-\rhomax)\juerg{^+}|^2
  \non
  \\
  && {} + \intQt F'(\rho) \, (\rho-\rhomax)^+
  + \intSt \FG'(\rhoG) \, (\rhoG-\rhomax)^+
  = \intQt \mu g'(\rho) (\rho-\rhomax)^+ \,.
  \non
\Esist
All \juerg{of} the terms on the \lhs\ are nonnegative, \juerg{while} the \rhs\ is nonpositive.
We conclude that $(\rho(t)-\rhomax)^+=0$ \juerg{\,in $\,\overline\Omega\,$} for every $t>0$, i.e., our assertion.

\step
Consequence

Since $g$, $F$ and $\FG$ are smooth on $(-1,1)$ 
and \eqref{hpg} implies that $g$ is strictly positive on~$(-1,1)$,
the separation inequalities \eqref{separaz} imply the bounds
\Beq
  g(\rho) \geq \gstar > 0 
  \aand
  |\Phi(\rho)| \leq C^* 
  \quad \hbox{in \gianni{$\overline Q_\infty$}},
  \quad
  |\Phi_\Gamma(\rhoG)| \leq C^*
  \quad \hbox{on \gianni{$\overline\Sigma_\infty$}},
  \label{daseparaz}
\Eeq
for $\Phi\in\{g,g'\!,g''\!,F,F'\!,F''\}$ and $\Phi_\Gamma\in\{\FG,\FG',\FG''\}$,
and for some constants $\gstar$ and $C^*$ that \juerg{depend only}
on the shape of the nonlinearities and the initial datum~$\rhoz$.
In particular, they do not depend on~$\eps$.

\step 
Uniqueness

We prove that the solution to problem \Pbl\ with $\eps=0$ is unique.
To this end, we fix $T>0$ and two solutions $(\mu_i,\rho_i,\rho_\juerg{{i_\Gamma}})$, $i=1,2$,
and show that they coincide on~$[0,T]$.
We set for convenience $\mu:=\mu_1-\mu_2$ and analogously define $\rho$ and~$\rhoG$.
Then, we write \eqref{prima} for both solutions
and test the difference by~$\mu$.
Using the identity
\Bsist
  && \{ 
    2g(\rho_1) \dt\mu_1 
    + \mu_1 g'(\rho_1) \dt\rho_1
    - 2g(\rho_2) \dt\mu_2
    - \mu_2 g'(\rho_2) \dt\rho_2
  \} \mu
  \non
  \\
  && = \dt \bigl( g(\rho_1) \, \mu^2 \bigr)
  + 2 \dt\mu_2 \, \bigl( g(\rho_1) - g(\rho_2) \bigr) \, \mu
  + \mu_2 \bigl( g'(\rho_1) \dt\rho_1 - g'(\rho_2) \dt\rho_2 \bigr) \mu\,, 
  \non
\Esist
we obtain \juerg{that}
\Bsist
  && \iO \juerg{g(\rho_1(t))} \, |\mu(t)|^2
  + \intQt |\nabla\mu|^2 
  \non
  \\
  && = - \intQt 2 \dt\mu_2 \, \bigl( g(\rho_1) - g(\rho_2) \bigr) \, \mu
  - \intQt \mu_2 \bigl( g'(\rho_1) \dt\rho_1 - g'(\rho_2) \dt\rho_2 \bigr) \mu \,.
  \label{testdiffprima}
\Esist
Next, we write \eqref{seconda} at the time~$s$ for both solutions,
test the difference by~$\dt(\rho,\rhoG)(s)$, and integrate over~$(0,t)$ with respect to~$s$.
Then, we add $\intQt\rho\,\dt\rho+\intSt\rhoG\,\dt\rhoG$ to both sides.
We~get
\Bsist
  && \intQt |\dt\rho|^2
  + \intSt |\dt\rhoG|^2
  + \frac 12 \, \normaV{\rho(t)}^2
  + \frac 12 \, \normaVG{\rhoG(t)}^2
  \non
  \\
  && = - \intQt \bigl(F'(\rho_1) - F'(\rho_2) \bigr) \dt\rho
  - \intSt \bigl( \juerg{F'_\Gamma(\rho_{1_\Gamma}) - F'_\Gamma(\rho_{2_\Gamma})} \bigr) \dt\rhoG
  \non
  \\
  && \quad {}
  + \intQt \bigl( \mu_1 g'(\rho_1) - \mu_2 g'(\rho_2) \bigr) \dt\rho
  + \intQt \rho \, \dt\rho 
  + \intSt \rhoG \, \dt\rhoG \,.
  \label{testdiffseconda}
\Esist
At this point, we add \accorpa{testdiffprima}{testdiffseconda} to each other
and use the separation property, the first \pier{inequality in} \eqref{daseparaz} for~$\rho_1$,
and the boundedness and the \Lip\ continuity of the nonlinearities on~$[\rhomin,\rhomax]$.
We~\juerg{find that}
\Bsist
  && \gstar \! \iO |\mu(t)|^2
  + \intQt |\nabla\mu|^2 
  + \intQt |\dt\rho|^2
  + \intSt |\dt\rhoG|^2
  + \frac 12 \, \normaV{\rho(t)}^2
  + \frac 12 \, \normaVG{\rhoG(t)}^2
  \non
  \\
  \separa
  && \leq c \intQt |\dt\mu_2| \, |\rho| \, |\mu|
  + c \intQt \mu_2 \bigl( |\dt\rho| + |\rho| \, |\dt\rho_2| \bigr) \, |\mu|
  \non
  \\
  && \quad {}
  + c \intQt |\rho| \, |\dt\rho|
  + c \intSt |\rhoG| \, |\dt\rhoG|
  + c \intQt \bigl( \mu_1 |\rho| + |\mu| \bigr) \, |\dt\rho| \,.
  \label{perunicita}
\Esist
Many integrals on the \rhs\ can be dealt with just using the \Holder\ and Young inequalities.
Thus, we consider just the terms that need some treatment.
In the next lines, we owe to the continuous embeddings $V\subset\Lx p$ for $p\in[1,6]$ and $W\subset\Cx0$,
and $\delta$~is a positive parameter.
We~have
\Bsist
  && \intQt |\dt\mu_2| \, |\rho| \, |\mu|
  \leq \iot \norma{\dt\mu_2(s)}_2 \norma{\rho(s)}_4 \norma{\mu(s)}_4 \, ds
  \non
  \\
  && \leq \delta \iot \normaV{\mu(s)}^2 \, ds
  + c_\delta \iot \normaH{\dt\mu_2(s)}^2 \normaV{\rho(s)}^2 \, ds\,, 
  \non
\Esist
and we notice that the function $s\mapsto\normaH{\dt\mu_2(s)}^2$
belongs to $L^1(0,T)$ by \eqref{regmu} for~$\mu_2$.
We estimate the next integral as follows,
\Bsist
  && \intQt \mu_2 \bigl( |\dt\rho| + |\rho| \, |\dt\rho_2| \bigr) \, |\mu|
  \non
  \\
  && \leq \iot \norma{\mu_2(s)}_\infty \norma{\dt\rho(s)}_2 \norma{\mu(s)}_2 \, ds
  + \juerg{c}\iot \norma{\mu_2(s)}_6 \norma{\rho(s)}_6 \norma{\dt\rho_2(s)}_6 \norma{\mu(s)}_6 \, ds
  \non
  \\
  && \leq \delta \iot \normaH{\dt\rho(s)}^2 \, ds
  + c_\delta \iot\normaW{\mu_2(s)}^2 \normaH{\mu(s)}^2 \, ds 
  \non
  \\
  && \quad {}
  + \delta \iot \normaV{\mu(s)}^2 \, ds
  + c_\delta \iot\normaV{\mu_2(s)}^2 \normaV{\dt\rho_2(s)}^2 \normaV{\rho(s)}^2 \, ds\,, 
  \non
\Esist
and we point out that the functions 
$s\mapsto\normaW{\mu_2(s)}^2$, $s\mapsto\normaV{\mu_2(s)}^2$, and $s\mapsto\normaV{\dt\rho_2(s)}^2$,
belong to $L^1(0,T)$, $L^\infty(0,T)$ and $L^1(0,T)$, respectively,
due to \accorpa{regmu}{regrho} for $\mu_2$ and~$\rho_2$.
Finally, we estimate one \juerg{further} term.
We~have \juerg{that}
\Bsist
  && \intQt \mu_1 |\rho| \, |\dt\rho| 
  \leq \iot \norma{\mu_1(s)}_4 \norma{\rho(s)}_4 \norma{\dt\rho(s)}_2 \, ds
  \non
  \\
  && \leq \delta \iot \normaH{\dt\rho}^2 \, ds
  + c_\delta \iot \normaV{\mu_1(s)}^2 \normaV{\rho(s)}^2 \, ds\,,
  \non
\Esist
\juerg{where} the function $s\mapsto\normaV{\mu_1(s)}^2$ belongs to $L^\infty(0,T)$.
Therefore, by choosing $\delta$ small enough and coming back to~\eqref{perunicita},
we can apply the Gronwall lemma to conclude that $(\mu,\rho,\rhoG)$ vanishes on
\juerg{$\overline\Omega\times[0,T]$}.

\medskip

Now, we show the existence of a solution to problem \Pbl\ with $\eps=0$
and prove the last sentence of the statement of Theorem~\ref{Wellposedness}.
To do that, it suffices to establish a number of a~priori estimates 
on the solution $(\mueps,\rhoeps,\rhoGeps)$ on an arbitrarily fixed time interval $[0,T]$
and to use proper compactness results.
\juerg{As the} uniqueness of the solution to the limiting problem \juerg{is} already known,
\juerg{it follows that the convergence properties proved below for a subsequence actually hold} for the whole family.
In view of the asymptotic \bhv\ that we aim to study in the next section,
we distinguish in the notation the constants that may depend on~$T$,
as explained at the end of Section~\ref{STATEMENT}.
Of course, we can assume $\eps\leq 1$.
In order to keep the length of the paper reasonable,
we perform some of the next estimates just formally.

\step
First a priori estimate

We observe that
\Beq
  \bigl\{ (\eps + 2g(\rhoeps))\dt\mueps + \mueps g'(\rhoeps)\dt\rhoeps \bigr\} \mueps
  = \dt\bigl(\bigl( {\textstyle\frac\eps 2} + g(\rhoeps) \bigr) |\mueps|^2 \bigr).
  \non
\Eeq
Hence, if we multiply \eqref{prima} by $\mueps$ and integrate over~$Q_t$, we obtain that
\Beq
  \frac \eps 2 \iO |\mueps(t)|^2 
  + \iO g(\rhoeps(t)) |\mueps(t)|^2
  + \intQt |\nabla\mueps|^2
  = \frac \eps 2 \iO \muz^2 + \iO g(\rhoz) \muz^2 \,.
  \non
\Eeq
By accounting for \pier{\eqref{separaz} and} \eqref{daseparaz}, we deduce, 
for every $t\geq0$, the global estimate 
\Beq
  \gstar \! \iO |\mueps(t)|^2
  + \intQt |\nabla\mueps|^2
  \leq \frac 12 \iO \muz^2 + \iO g(\rhoz) \muz^2 
  = c \,.
  \label{primastima}
\Eeq

\step
Second a priori estimate

We write \eqref{seconda} at the time~$s$
and \pier{choose the test pair $(v,v_\Gamma) = (\dt\rhoeps,\dt\rhoGeps)(s)$,
which is allowed by the regularity \eqref{regrho}.}
Then, we integrate over~$(0,t)$.
\pier{Thanks} to the Schwarz and Young inequalities, we~have
\Bsist
  && \intQt |\dt\rhoeps|^2 
  + \intSt |\dt\rhoGeps|^2
  + \frac 12 \iO |\nabla\rhoeps(t)|^2
  + \frac 12 \iG |\nablaG\rhoGeps(t)|^2
  \non
  \\
  && \quad {}
  + \iO F(\rhoeps(t))
  + \iG \FG(\rhoGeps(t))
  \non
  \\
  && = \frac 12 \iO |\nabla\rhoz|^2
  + \frac 12 \iG |\nablaG\rhozG|^2
  + \iO F(\rhoz)
  + \iG \FG(\rhozG)
  + \intQt \mueps g'(\rhoeps) \dt\rhoeps
  \non
  \\
  && \leq c
  + \frac 12 \intQt |\dt\rhoeps|^2
  + c \intQt |\mueps|^2 .
  \non
\Esist
Since $|\rhoeps|\leq 1$, \eqref{primastima}~holds, 
and \eqref{hpFsecondo} implies that $F$ and $\FG$ are bounded from below, 
we deduce that
\Beq
  \norma{(\rhoeps,\rhoGeps)}_{\H1\calH\cap\L\infty\calV}
  + \norma{F(\rhoeps)}_{\L\infty\Luno}
  + \norma{\FG(\rhoGeps)}_{\L\infty\LunoG}
  \leq c_T \,.
  \label{secondastima}
\Eeq

\step
Third a priori estimate

\gianni{By starting from \Secondabvp\ 
and accounting for \eqref{daseparaz} and \accorpa{primastima}{secondastima},
we successively deduce a number of estimates
with the help of the inequalities \Richiami, written with $v=\rhoeps(t)$ and 
$\vG=\rhoGeps(t)$ \juerg{and 
then squared and integrated} over~$(0,T)$.
We have
\Bsist
  && \norma{\Delta\rhoeps}_{\L2H} \leq c_T 
  \quad \hbox{from \eqref{secondaO}},
  \non
  \\
  && \norma{\dn\rhoeps}_{\L2{\HxG{-1/2}}} \leq c_T 
  \quad \hbox{from \eqref{dnHmum}},
  \non
  \\
  \separa
  && \norma{\DeltaG\rhoGeps}_{\L2{\HxG{-1/2}}} \leq c_T 
  \quad \hbox{from \eqref{secondaG}},
  \non
  \\
  && \norma\rhoGeps_{\L2{\HxG{3/2}}} \leq c_T  
  \quad \hbox{from \eqref{stimaHtmG}},
  \non
  \\
  \separa
  && \norma\rhoeps_{\L2\Hdue} \leq c_T 
  \quad \hbox{from \eqref{stimaHdueD}},
  \non
  \\
  && \norma{\dn\rhoeps}_{\L2\HG} \leq c_T  
  \quad \hbox{from \eqref{dnLd}},
  \non
  \\
  \separa
  && \norma{\DeltaG\rhoGeps}_{\L2\HG} \leq c_T 
  \quad \hbox{from \eqref{secondaG}},
  \non
  \\
  && \norma\rhoGeps_{\L2\HdueG} \leq c_T 
  \quad \hbox{from \eqref{stimaHdueG}}.
  \non
\Esist
In conclusion, we have proved that}
\Beq
  \norma\rhoeps_{\L2\Hdue} 
  + \norma\rhoGeps_{\L2\HdueG} \leq c_T \,.
  \label{quartastima}
\Eeq

\step
\juerg{Fourth a priori estimate}

We (formally) differentiate \eqref{seconda} with respect to time
and set $\zeta:=\dt\rhoeps$ and $\zetaG:=\dt\rhoGeps$, for brevity.
Then we write the variational equation we obtain at the time~$s$
and test it by $(\zeta,\zetaG)(s)$.
Finally, we integrate over $(0,t)$ and add $C\intQt|\zeta|^2+C\intSt|\zetaG|^2$ to both sides,
where $C$ is the constant that appears in~\eqref{hpFsecondo}.
We obtain \juerg{the identity}
\Bsist
  && \frac 12 \iO |\zeta(t)|^2
  + \frac 12 \iG |\zetaG(t)|^2
  + \intQt |\nabla\zeta|^2
  + \intSt |\nablaG\zetaG|^2
  \non
  \\
  && \quad {}
  + \intQt \bigl( F''(\rhoeps) + C \bigr) |\zeta|^2
  + \intQt \bigl( \FG''(\rhoGeps) + C \bigr) |\zetaG|^2
  \non
  \\
  && = \frac 12 \iO |\zeta(0)|^2
  + \frac 12 \iG |\zetaG(0)|^2
  + \intQt \dt\mueps \, g'(\rhoeps) \zeta
  + \intQt \mueps g''(\rhoeps) |\zeta|^2 
  \non
  \\
  && \quad {}
  + C \intQt |\zeta|^2
  + C \intSt |\zetaG|^2 \,.
  \label{perquintastima}
\Esist
All \juerg{of} the terms on the \lhs\ are nonnegative, \juerg{while the second volume integral over $\,Q_t\,$
on the \rhs\ is nonpositive since $\mueps\geq0$ and $g''\leq0$.}
It remains to find bounds for \juerg{the first volume integral
over $\,Q_t\,$ on the \rhs\ and for the sum of the terms that involve the initial values. We handle the latter first.} 
To this end, we write \eqref{seconda} at the time $t=0$ and test it by 
$(v,\vG)=(\zeta,\zetaG)(0)$.
We obtain
\begin{eqnarray}
  && \iO |\zeta(0)|^2
  + \iG |\zetaG(0)|^2
  = -\iO \nabla\rhoz \cdot \nabla\zeta(0)
  - \iG \nablaG\rhozG \cdot \nablaG\zetaG(0)
  \non
  \\
	  &&\label{juerg1} \quad {}
  - \iO F'(\rhoz) \zeta(0) 
  - \iG \FG'(\rhozG) \zetaG(0)
  + \iO \muz g'(\rhoz) \zeta(0) \,.
  \end{eqnarray}
On account of \eqref{hprhoz}, we have, \juerg{using Young's inequality \pier{and \eqref{dnLd}}},
\Bsist
  && -\iO \nabla\rhoz \cdot \nabla\zeta(0)
  - \iG \nablaG\rhozG \cdot \nablaG\zetaG(0)
  =  \iO \Delta\rhoz \, \zeta(0)
  - \iG \bigl( \dn\rhoz - \DeltaG\rhozG \bigr) \zetaG(0)
  \non
  \\
  && \leq \juerg{\frac 14} \iO |\zeta(0)|^2
  + \juerg{\frac 14} \iG |\zetaG(0)|^2
  + c \,\norma\rhoz_{\Hx2}^2
  + c\, \norma\rhozG_{\HxG2}^2 \,.
    \non
\Esist
\juerg{Moreover, it follows from \eqref{hpmuz}, \eqref{hprhoz}, \eqref{daseparaz},
and Young's inequality that the expression in the second line of \eqref{juerg1}
is bounded by} 
\begin{equation*}
\juerg{\frac 14\iO|\zeta(0)|^2\,+\,\frac 14\iG|\zetaG(0)|^2\,+\,c\,.}
\end{equation*}
\juerg{We thus have shown that}
\begin{equation}
\label{juerg2}
\juerg{\iO |\zeta(0)|^2   + \iG |\zetaG(0)|^2\,\le\,c\,.}
\end{equation}
\juerg{It remains to bound the first volume integral over $Q_t$ in 
\eqref{perquintastima}, which we denote by $I$. 
This estimate requires more effort. 
At first, observe that \gianni{\eqref{prima}} implies that
\begin{equation}
\label{juerg3}
\dt\mueps\,=\,\frac 1{\eps+2g(\rhoeps)}\,\Delta\mueps\,-\,
\frac{g'(\rhoeps)}{\eps+2g(\rhoeps)}\,\zeta\,\mueps\,,
\end{equation}
where, thanks to \eqref{daseparaz}, $1/(\eps+2g(\rhoeps))\le 1/(2g^*)$ for all $\eps>0$. 
Now, using \eqref{juerg3}, we find that
\begin{equation}
\label{juerg4}
I\,=\intQt \frac{g'(\rhoeps)\,\zeta}{\eps+2g(\rhoeps)}\,\Delta\mueps\,-\,
\intQt\mueps\,\frac{(g'(\rhoeps))^2}{\eps+2g(\rhoeps)}\,\zeta^2\,=:\,I_1+I_2\,,
\end{equation}
}%
\juerg{%
with obvious notation. 
The second integral is easy to handle. 
In fact, thanks to \eqref{daseparaz}, \eqref{compact}, and H\"older's and Young's inequalities,
we \pier{infer} that
\begin{eqnarray}
I_2&\!\!\le\!\!&c\int_0^t\|\mueps(s)\|_4\,\|\zeta(s)\|_2\,\|\zeta(s)\|_4\,ds\non\\
&\!\!\le\!\!& \frac 16 \intQt|\nabla\zeta|^2\,+\,c\int_0^t
\left(1+\|\mueps(s)\|_V^2\right)\|\zeta(s)\|_H^2\,ds\,,\label{juerg5}
\end{eqnarray}
}%
\juerg{where we know from \eqref{primastima} that $\,\int_0^T\|\mueps(s)\|_V^2\,ds\le c_T$. 
For the first integral, integration by parts and \eqref{daseparaz} yield that}
\begin{eqnarray}
\juerg{I_1}&\juerg{\!\!=\!\!}&\juerg{-\intQt\nabla\mueps\cdot\nabla\Big(\frac{g'(\rhoeps)\,\zeta}
{\eps+2g(\rhoeps)}\,\Big)}\non\\
&\!\!\le\!\!&\juerg{ C_1\intQt|\nabla\mueps|\,|\nabla\zeta|\,+\,C_1\intQt|\nabla\mueps|\,
|\nabla\rhoeps|\,|\zeta|\,=:\,C_1(I_{11}+I_{12}),}\label{juerg6}
\end{eqnarray}
\juerg{with obvious notation. 
Clearly, owing to \eqref{primastima} and Young's inequality,} \pier{we find that}
\begin{equation}
\label{juerg7}
\juerg{C_1\,I_{11}\,\le\,\frac 16\intQt|\nabla\zeta|^2\,+\,c\,.}
\end{equation}
\juerg{Moreover, invoking H\"older's and Young's inequalities, 
\gianni{the compactness inequality \eqref{compact},
as well as the continuity of the embedding $H^2(\Omega)\subset W^{1,4}(\Omega)$},  
we \pier{infer} that}
\begin{eqnarray}
C_1\,I_{12}&\!\!\le\!\!&C_1\int_0^t\|\nabla\mueps(s)\|_2\,\|\nabla\rhoeps(s)\|_4
\,\|\zeta(s)\|_4\,ds\non\\
&\!\!\le\!\!&\frac 16\intQt|\nabla\zeta|^2\,+\,c\intQt|\zeta|^2\,+\,
c\int_0^t\|\nabla\mueps(s)\|_2^2\,\|\rhoeps(s)\|^2_{\Hdue}\,ds\,.
\label{juerg8}
\end{eqnarray}
\juerg{Notice that $\,\int_0^T\|\nabla\mueps(s)\|_2^2\,ds\le c\,$ for every $T>0$,
by virtue of \eqref{primastima}. 
We now aim to estimate $\|\rhoeps(s)\|_{\Hdue}$ in terms of $\zeta$ and $\zetaG$. 
To this end, we derive a chain of estimates which are each valid for almost every $s\in (0,T)$. 
To begin with,
we deduce from \eqref{primastima} and \eqref{secondastima} that}
\begin{equation}
\label{juerg9}
\juerg{\|\Delta\rhoeps(s)\|_2\,=\,\|\zeta(s)+F'(\rhoeps(s))-\mueps(s)\,g'(\rhoeps(s))\|_2
\,\le\,c\,+\,\|\zeta(s)\|_2\,.}
\end{equation}
\juerg{Consequently, by \eqref{dnHmum} \pier{we have that}}
\begin{equation}
\label{juerg10}
\juerg{\|\dn\rhoeps(s)\|_{H^{-1/2}(\Gamma)}\,\le\,C_\Omega\left(\|\rhoeps(s)\|_V+\|
\Delta\rhoeps(s)\|_2\right)\le\,c_T\,(1+\|\zeta(s)\|_2)\,,}
\end{equation}
\juerg{and \gianni{\eqref{secondaG}}, \eqref{daseparaz} and \eqref{secondastima} imply that
\begin{eqnarray}
\|\Delta_\Gamma\rhoGeps(s)\|_{H^{-1/2}(\Gamma)}&\!\!\le\!\!&\|\dn\rhoeps(s)+F'_\Gamma
(\rhoGeps(s))+\zeta_\Gamma(s)\|_{H^{-1/2}(\Gamma)}\non\\[1mm]
&\!\!\le\!\!&
c_T\,(1+\|\zeta(s)\|_2+\|\zetaG(s)\|_{\HG})\,.
\label{juerg11}
\end{eqnarray} 
}
\juerg{
But then, thanks to \eqref{stimaHtmG} and \eqref{secondastima} \pier{it is clear that}
\begin{eqnarray}
\|\rhoGeps(s)\|_{H^{3/2}(\Gamma)}&\!\! \leq\!\!&\CO \left( 
\|\rhoGeps(s)\|_{\HunoG} + \|\DeltaG\rhoGeps(s)\|_{H^{-1/2}(\Gamma)} \right)
  \non\\[1mm]
	&\!\!\le\!\!&c_T\,(1+\|\zeta(s)\|_2+\|\zetaG(s)\|_{\HG})\,,
\end{eqnarray}
}
\juerg{	
whence, owing to \eqref{stimaHdueD}, we finally arrive at the estimate
\begin{equation}
\label{juerg13}
\|\rhoeps(s)\|_{\Hdue}\,\le\,c_T\left(1+\|\zeta(s)\|_H+\|\zetaG(s)\|_{\HG}\right)\,.
\end{equation}
}
\juerg{We thus obtain from \eqref{juerg8} that
\begin{equation}
C_1\,I_{12}\,\le\,\frac 16\intQt|\nabla\zeta|^2\,+\,c\intQt|\zeta|^2\,+\,c_T
\,+\,c_T\int_0^t\|\nabla\mueps(s)\|_2^2\left(\|\zeta(s)\|_H^2+\|\zetaG(s)\|^2
_{\HG}\right) ds\,.
\label{juerg14}
\end{equation}
}%
\juerg{Therefore, recalling \eqref{perquintastima} and invoking the estimates
\eqref{juerg2}, \eqref{juerg5}--\eqref{juerg8}, we can apply
Gronwall's lemma and conclude that}
\Beq
  \norma{(\dt\rhoeps,\dt\rhoGeps)}_{\L\infty\calH\cap\L2\calV}
  \leq c_T \,.
  \label{quintastima}
\Eeq

\step
Fifth a priori estimate

\juerg{We now notice that \eqref{juerg13} and \eqref{quintastima} imply that}
\begin{equation}
\label{juerg15}
\juerg{\|\rhoeps\|_{L^\infty(0,T;\Hdue)}\,\le\,c_T\,.}
\end{equation}
\juerg{Then we may infer from \eqref{dnLd}, \gianni{\eqref{secondaG}}, \eqref{stimaHdueG}, in this order, the estimates}
\begin{equation*}
\juerg{\|\dn\rhoeps\|_{L^\infty(0,T;\HG)}\le c_T, \quad \|\Delta_\Gamma\rhoGeps\|_{L^\infty(0,T;\HG)}\le c_T,
\quad \|\rhoGeps\|_{L^\infty(0,T;\HdueG)}\le c_T,}
\end{equation*}
\juerg{so that}

\Beq
  \norma{(\rhoeps,\rhoGeps)}_{\L\infty{\Hdue\times\HdueG}}
  \leq c_T \,.
  \label{sestastima}
\Eeq

\step
Sixth a priori estimate

At this point, we can multiply \eqref{prima} by $\dt\mueps$ and integrate over~$Q_t$.
Then, we add $\intQt\mueps\dt\mueps$ to both sides.
By owing to the \Holder, Sobolev and Young inequalities,
we obtain
\Bsist
  && \intQt \bigl( \eps + 2g(\rhoeps) \bigr) |\dt\mueps|^2
  + \frac 12 \, \normaV{\mueps(t)}^2
  \non
  \\
  && = \frac 12 \, \normaV\muz^2
  + \intQt \mueps \dt\mueps
  - \intQt \mueps g'(\rhoeps) \dt\rhoeps \dt\mueps 
  \non
  \\
  && \leq c + \iot \norma{\mueps(s)}_2 \norma{\dt\mueps(s)}_2 \, ds
  + c \iot \norma{\mueps(s)}_4 \norma{\dt\rhoeps(s)}_4 \norma{\dt\mueps(s)}_2 \, ds
  \non
  \\
  && \leq c + \gstar \intQt |\dt\mueps|^2
  + c \, \norma\mueps_{L^2(0,t;H)}^2
  + c \iot \normaV{\dt\rhoeps(s)}^2 \normaV{\mueps(s)}^2 \, ds\,, 
  \non
\Esist
where $\gstar$ is the constant \juerg{introduced} in~\eqref{daseparaz}.
As $\,2g(\rhoeps)\geq2\gstar$, \juerg{we  may use \eqref{primastima}, \eqref{quintastima}
and Gronwall's lemma to conclude~that}
\Beq
  \norma\mueps_{\H1H\cap\L\infty V} \leq c_T \,.
  \label{settimastima}
\Eeq
By comparison in \eqref{prima}, we estimate $\Delta\mueps$.
Hence, by applying~\eqref{stimaHdueN}, we derive~that
\Beq
  \norma\mueps_{\L2W} \leq c_T \,.
  \label{dasettimastima}
\Eeq

\step
Conclusion

If we collect all the previous estimates and use standard compactness results,
then we have (\pier{in principle} for a subsequence) \juerg{that}
\Bsist
  & \mueps \to \mu
  & \quad \hbox{in $\H1H\cap\L\infty V\cap\L2W$}\,,
  \non
  \\
  & (\rhoeps,\rhoGeps) \to (\rho,\rhoG)
  & \quad \hbox{in $\W{1,\infty}\calH\cap\H1\calV\cap\L\infty{\Hdue\times\HdueG}$}\,,
  \non
\Esist
as $\eps\searrow 0$,
the convergence being understood in the sense of the corresponding weak star topologies.
Notice that the limiting triplet fulfills the regularity requirements \Regsoluz.
Next, by the compact embeddings $V\subset\Lx5$, $\Hdue\subset\Cx0$, and $\HdueG\subset\CxG0$,
and using well-known strong compactness results
(see, e.g., \cite[Sect.~8, Cor.~4]{Simon}),
we deduce the useful strong convergence
\Beq
  \mueps \to \mu
  \quad \hbox{in $\C0{\Lx5}$}
  \aand
  (\rhoeps,\rhoGeps) \to (\rho,\rhoG)
  \quad \hbox{in \juerg{$C^0(\overline Q)\times C^0(\overline \Sigma)$}}.
  \label{pier1}
\Eeq
This allows us to deal with nonlinearities
and to take the limits of the products that appear in the equations.
Hence, we easily conclude that the triplet $(\mu,\rho,\rhoG)$
solves \pier{\eqref{prima} and the time-integrated version of \eqref{seconda} on $(0,T)$
(which is equivalent to \eqref{seconda} itself) with $\eps=0$. Moreover, the initial conditions \eqref{cauchy} easily pass to the limit in view of \eqref{pier1}.}
This concludes the existence proof.
By uniqueness, the whole family $\{(\mueps,\rhoeps,\rhoGeps)\}$ converges to $(\mu,\rho,\rhoG)$
in the above topology as $\eps\searrow 0$.\QED


\section{\Longtime\ \bhv}
\label{LONGTIME}
\setcounter{equation}{0}

This section is devoted to the proof of Theorem~\ref{Omegalimit}.
In the sequel, it is understood that $\eps\in[0,1]$ is fixed
and that $(\mu,\rho,\rhoG)$ is the unique solution to problem \Pbl\ 
given by Theorems~\ref{Wellposednesseps} and~\ref{Wellposedness}
in the two cases $\eps>0$ and $\eps=0$, respectively.
First of all, we have to show that the \omegalimit~\eqref{omegalimit} is \juerg{nonempty}.
This \juerg{necessitates} proper a~priori estimates on the whole half-line $\{t\geq0\}$.

\step
First global estimate

From~\eqref{primastima}, we immediately deduce that
\Beq
  \norma\mu_{\LL\infty H} \leq c
  \aand
  \intQi |\nabla\mu|^2 \leq c \,.
  \label{primaglob}
\Eeq

\step
Second global estimate

We start by rearranging \eqref{prima} as follows:
\Beq
  \mu g'(\rho) \dt\rho
  = \dt \bigl( (\eps+2g(\rho)) \mu \bigr)
  - \Delta\mu \,.
  \label{perprimaid}
\Eeq
Now, we test \eqref{seconda}, written at the time~$s$, by $\dt(\rho,\rhoG)(s)$,
integrate over~$(0,t)$ and replace the \rhs\ \pier{with the help of}~\eqref{perprimaid}.
We~\juerg{obtain the identity}
\Bsist
  && \intQt |\dt\rho|^2
  + \intSt |\dt\rhoG|^2
  + \frac 12 \iO |\nabla\rho(t)|^2
  + \frac 12 \iG |\nablaG\rhoG(t)|^2
  \non
  \\
  && \quad {}
  + \iO F(\rho(t))
  + \iG \FG(\rhoG(t))
  \non
  \\
  && = \frac 12 \iO |\nabla\rhoz|^2
  + \frac 12 \iG |\nablaG\rhozG|^2
  + \iO F(\rhoz)
  + \iG \FG(\rhozG)
  + \intQt \mu g'(\rho) \dt\rho
  \non
  \\
  && = c + \iO \bigl( \eps + 2g(\rho(t)) \bigr) \mu(t)
  - \iO \bigl( \eps + 2g(\rhoz) \bigr) \muz
  - \intQt \Delta\mu \,.
  \non
\Esist
The last integral vanishes since $\dn\mu=0$.
By recalling that $F$ and $\FG$ are bounded from below and that $|\rho|\leq 1$,
and using~\eqref{primaglob}, we deduce that
\Beq
  \norma{(\rho,\rhoG)}_{\LL\infty\calV} \leq c \,, \quad
  \intQi |\dt\rho|^2 \leq c
  \aand
  \intSi |\dt\rhoG|^2 \leq c \,.
  \label{secondaglob}
\Eeq

\step
First conclusion

The first inequalities of \eqref{primaglob} and \eqref{secondaglob}\pier{, along
with the continuity of $(\mu,\rho,\rhoG)$ from $[0,+\infty)$ to $H\times\calV$,}
ensure that the \omegalimit~\eqref{omegalimit} is \juerg{nonempty}.
Namely, every divergent sequence of times \juerg{contains} a subsequence $\tn\nearrow+\infty$
such that $(\mu,\rho,\rhoG)(\tn)$ converges weakly in~$H\times\calV$.

\juerg{After establishing the first part of Theorem~\ref{Omegalimit},
we prove the second one}.
Thus, we pick any element $(\muo,\rhoo,\rhoGo)$ of the \omegalimit~\eqref{omegalimit}
and show that it is a stationary solution 
of our problem, i.e., that $\muo$ is a constant $\mus$ and that
the pair $(\rhoo,\rhoGo)$ coincides with a solution
$(\rhos,\rhoGs)$ to problem~\eqref{secondastaz}.
To this end, we fix a sequence $\tn\nearrow+\infty$ such~that
\Beq
  (\mu,\rho,\rhoG)(\tn) \to (\muo,\rhoo,\rhoGo)
  \quad \hbox{weakly in $H\times\calV$}
  \label{toomega}
\Eeq
and study the \bhv\ of the solution on the time interval $[\tn,\tn+T]$ 
with a fixed $T>0$.
For convenience, we shift everything to~$[0,T]$ 
by introducing $(\mun,\rhon,\rhoGn):[0,T]\to H\times\calV$ as follows
\Beq
  \mun(t) := \mu(\tn+t) , \quad
  \rhon(t) := \rho(\tn+t)
  \aand
  \rhoGn(t) := \rhoG(\tn+t)
  \quad \hbox{for $t\in[0,T]$} \,.
  \label{soluzn}
\Eeq
As $T$ is fixed once and for all, we do not care 
on the dependence of the constants on $T$ even in the notation,
and write $Q$ and $\Sigma$ for $Q_T$ and~$\Sigma_T$, respectively.
The inequalities \eqref{primaglob} and \eqref{secondaglob} imply that
\Bsist
  && \norma{(\mun,\rhon,\rhoGn)}_{\L\infty{H\times\calV}} \leq c\,,
  \label{perprimaconv}
  \\
  && \limn \Bigl( \intQ |\nabla\mun|^2 + \intQ |\dt\rhon|^2 + \intS |\dt\rhoGn|^2 \Bigr) = 0 \,.
  \label{perstaz}
\Esist
The bound \eqref{perprimaconv} yields a convergent subsequence in the weak star topology.
If we still label it by the index~$n$ to simplify the notation, we have
\Beq
  (\mun,\rhon,\rhoGn) \to (\mui,\rhoi,\rhoGi)
  \quad \hbox{weakly star in $\L\infty{H\times\calV}$} .
  \label{primaconv}
\Eeq
Now, we aim to improve the quality of the convergence.
Thus, we derive further estimates.

\step
First auxiliary estimate

A~partial use of \eqref{perstaz} provides a bound, namely
\Beq
  \norma\mun_{\L2V}
  + \norma{(\dt\rhon,\dt\rhoGn)}_{\L2\calH}
  \leq c \,.
  \label{primaaux}
\Eeq

\step
Second auxiliary estimate

We can repeat the argument that led to \eqref{quartastima}
and arrive~at
\Beq
  \norma{(\rhon,\rhoGn)}_{\L2{\Hdue\times\HdueG}} \leq c \,.
  \label{secondaaux}
\Eeq

\step
Third auxiliary estimate

We recall that $\mun$
and the space derivatives $D_i\rhon$ and $D_ig(\rhon)=\juerg{g'(\rhon)}D_i\rhon$
are bounded in $\L\infty H\cap\L2{\Lx6}$, by
\eqref{perprimaconv}, \eqref{primaaux}, \eqref{secondaaux}, and the continuous embedding $V\subset\Lx6$.
On the other hand, the continuous embedding
\Beq
  \L\infty H \cap \L2{\Lx6}
  \subset 
  \L4{\Lx3} \cap \L6{\Lx{18/7}}
  \non
\Eeq
holds \juerg{true, by virtue of} the interpolation inequalities.
Therefore, we conclude~that
\Beq
  \norma\mun_{\L4{\Lx3}}
  + \norma{\nabla\rhon}_{\L4{\Lx3}}
  + \norma{\nabla g(\rhon)}_{\L6{\Lx{18/7}}}
  \leq c \,.
  \label{terzaaux}
\Eeq

\step
Fourth auxiliary estimate

We want to improve the convergence of~$\mun$.
However, we cannot multiply \eqref{prima} by $\dt\mu$ 
since we do not have any information on~$\nabla\mu(\tn)$.
Therefore, we derive an estimate for $\dt\mun$ in a dual space.
By recalling that $g(\rho)\geq\gstar$ (see~\eqref{daseparaz}),
we divide both sides of \eqref{prima} by $\eps+2g(\rho)$,
Then, we take an arbitrary test function $v\in\L4V$, multiply the equality we obtain by~$v$,
integrate over~$\Omega\times(\tn,\tn+T)$ and rearrange.
We~get
\Beq
  \intQ \dt\mun \, v
  = - \intQ \frac {\mun g'(\rhon) \dt\rhon v} {\eps+2g(\rhon)}
  + \intQ \Delta\mun \, \frac v {\eps+2g(\rhon)}\,, 
  \non
\Eeq
and we now treat the \juerg{terms} on the \pier{\rhs\ separately}.
The first one is handled \juerg{using \Holder's inequality, namely,}
\Beq
  - \intQ \frac {\mun g'(\rhon) \dt\rhon v} {\eps+2g(\rhon)}
  \leq c \norma\mun_{\L4{\Lx3}} \norma{\dt\rhon}_{\L2\Ldue} \norma v_{\L4{\Lx6}} \,.
  \non
\Eeq
We integrate  the other term by parts and use the \Holder, Sobolev and Young inequalities as follows:
\Bsist
  && \intQ \Delta\mun \, \frac v {\eps+2g(\rhon)} 
  = - \intQ \nabla\mun \cdot \frac {(\eps+2g(\rhon)) \nabla v - 2 v g'(\rhon) \nabla\rhon} {(\eps+2g(\rhon))^2} 
  \non
  \\
  && \leq c \norma{\nabla\mun}_{\L2H} \norma v_{\L2V}
  + c \ioT \norma{\nabla\mun(s)}_2 \norma{v(s)}_6 \norma{\nabla\rhon(s)}_3 \, ds
  \non
  \\
  && \leq c \norma{\nabla\mun}_{\L2H} \bigl( \norma v_{\L2V} + \norma{\nabla\rhon}_{\L4{\Lx3}} \norma v_{\L4{\Lx6}} \bigr) \,.
  \non
\Esist
Therefore, we have for every $v\in\L4V$
\Bsist
  && \intQ \dt\mun \, v
  \leq c \norma\mun_{\L4{\Lx3}} \norma{\dt\rhon}_{\L2H} \norma v_{\L4V} 
  \non
  \\
  && \quad {}
  + c \norma{\nabla\mun}_{\L2H} \pier{{}\bigl( 1+ \norma{\nabla\rhon}_{\L4{\Lx3}} \bigr) \norma v_{\L4V}}  \,.
  \non
\Esist
Hence, on account of \eqref{primaglob}, \eqref{secondaglob} and~\eqref{terzaaux}, we conclude that
\Beq
  \norma{\dt\mun}_{\L{4/3}\Vp} \leq c \,.
  \label{quartaaux}
\Eeq

\step
Conclusion

By recalling the estimates \accorpa{primaaux}{secondaaux} and~\eqref{quartaaux},
we see that the convergence \eqref{primaconv} can be improved as follows:
\Bsist
  & \mun \to \mui
  & \quad \hbox{in $\W{1,4/3}\Vp\cap\L\infty H\cap\L2V$} \,,
  \non
  \\
  & (\rhon,\rhoGn) \to (\rhoi,\rhoGi)
  & \quad \hbox{in $\H1\calH\cap\L\infty\calV\cap\L2{\Hdue\times\HdueG}$}\,, 
  \non
  \non
\Esist
all in the \juerg{sense} of the corresponding weak star topologies.
Now, we prove that the limiting triple $(\mui,\rhoi,\rhoGi)$
solves problem \accorpa{prima}{seconda},
the first equation being understood in a generalized sense.
By \cite[Sect.~8, Cor.~4]{Simon} and the compact embeddings $\Hdue\subset V\subset H \pier{{}\subset V^*}$
and $\HdueG\subset\VG\subset\HG$, 
we also have (for a not relabeled subsequence)\pier{%
\begin{align}
  & \mun \to \mui
  \quad \hbox{strongly in $\C0{V^*}\cap\L2H $ and a.e.\ in $Q$}
  \label{pier2}
  \\[0.1cm]
  & (\rhon,\rhoGn) \to (\rhoi,\rhoGi)
  \quad \hbox{strongly in $\C0\calH\cap\L2\calV$ and a.e.\ on $Q\times\Sigma$}
  \label{pier3}
  \\[0.1cm]
  & \nabla g(\rhon) = g'(\rhon) \nabla\rhon \to g'(\rhoi) \nabla\rhoi = \nabla g(\rhoi)
  \quad \aeQ \,.
  \label{pier4}
\end{align}
}%
It follows that $(F'(\rhon),\FG'(\rhoGn))$ converges to $(F'(\rhoi),\FG'(\rhoGi))$
in $\L\infty\calH$, just by \Lip\ continuity.
This allows us to conclude that $(\rhoi,\rhoGi)$ solves the time-integrated version of~\eqref{seconda},
thus equation \eqref{seconda} itself.
As for \eqref{prima}, we \pier{recall \eqref{terzaaux} and} notice that $4<6$ and $2<18/7$\pier{. Then, with the help of \eqref{pier4} and the Egorov theorem, we} deduce~that
\Bsist
  && \nabla g(\rhon) \to \nabla g(\rhoi)
  \quad \hbox{strongly in $(\L4\Ldue)^3$,\quad whence} 
  \non
  \\
  && g(\rhon) \to g(\rhoi)
  \quad \hbox{strongly in $\L4V$}.
  \non
\Esist
Therefore, if we assume that $v\in\L\infty{\Wx{1,\infty}}$, we have that
\Bsist
  && g(\rhon) v \to g(\rhoi) v
  \quad \hbox{strongly in $\L4V$,\quad whence}
  \non
  \\
  && \intQ \bigl( \eps + 2g(\rhon) \bigr) \dt\mun \, v
  \to {}_{\L{4/3}\Vp} \< \dt\mui , \bigl( \eps + 2g(\rhoi) \bigr) v >_{\L4V} \,.
  \non
\Esist
On the other hand, from the convergence almost everywhere, we also have
\Beq
  g'(\rhon) \to g'(\rhoi)
  \quad \hbox{strongly in $\L4{\Lx6}$},
  \non
\Eeq
since $g'(\rhon)$ is bounded in $\LQ\infty$.
Moreover, \eqref{terzaaux}~implies that
$\mun$ converges to $\mui$ weakly in $\L4{\Lx3}$\pier{. On the other hand,}
\eqref{perstaz} yields the strong convergence of $\dt\rhon$ to~$0$ in~$\L2H$
(by~the way, $0$~must coincide with~$\dt\rhoi$).
We deduce that
\Beq
  \mun g'(\rhon) \dt\rhon \to \mui g'(\rhoi) \dt\rhoi
  \quad \hbox{weakly in $\LQ1$} \,.
  \non
\Eeq
Therefore, we conclude that
\Bsist
  && {}_{\L{4/3}\Vp} \< \dt\mui , \bigl( \eps + 2g(\rhoi) \bigr) v >_{\L4V}
  \non
  \\
  && \quad {}
  + \intQ \mui g'(\rhoi) \dt\rhoi \, v
  + \intQ \nabla\mui \cdot \nabla v
  = 0
  \label{varinf}
\Esist
for every $v\in\L\infty{\Wx{1,\infty}}$.
On the other hand, we know that $\mui\in\L4{\Lx3}$ by \eqref{terzaaux} and 
that $\dt\rhoi\in\L2H$.
Since $g'$ is bounded and the continuous embedding $V\subset\Lx6$ implies $\Lx{6/5}\subset\Vp$,
we also~have that
\Beq
  \mui g'(\rhoi) \dt\rhoi \in \L{4/3}{\Lx{6/5}} \subset \L{4/3}\Vp \,.
  \non  
\Eeq
Hence, by a simple density argument, we see that the variational equation \eqref{varinf} also holds \juerg{true} for every $v\in\L4V$.
At this point, we observe that \eqref{perstaz} \juerg{implies} that
\Beq
  \nabla\mui = 0 \,, \quad
  \dt\rhoi = 0
  \aand
  \dt\rhoGi = 0 \,.
  \label{derivatenulle}
\Eeq
In particular, \eqref{varinf} reduces~to
\Beq
  {}_{\L{4/3}\Vp} \< \dt\mui , \bigl( \eps + 2g(\rhoi) \bigr) v >_{\L4V}
  = 0
  \quad \hbox{for every $v\in\L4V$} 
  \non
\Eeq
and we easily infer that $\dt\mui=0$.
Indeed, the inequality $g(\rhon)\geq\gstar$ for every~$n$
implies $g(\rhoi)\geq\gstar$.
Thus, every $\phi\in C^\infty_c(Q)$ 
can be written as $\phi=\pier{(\eps+2g(\rhoi))}v$ \pier{for} some $v\in\L4V$ 
since $\nabla\rhoi\in(\L4H)^3$ by~\eqref{terzaaux}.
Therefore, $\dt\mui$ actually vanishes and we conclude that $\mui$ takes a constant value~$\mus$.
From~\eqref{derivatenulle} we also deduce that $(\rhoi,\rhoGi)$ is a time-independent pair $(\rhos,\rhoGs)$,
so that \eqref{seconda} reduces to~\eqref{secondastaz}.
Finally, we show that $(\mus,\rhos,\rhoGs)=(\muo,\rhoo,\rhoGo)$.
Indeed, \pier{\eqref{pier2} and \eqref{pier3}} imply that
\Beq
  \juerg{(\mun,(\rhon,\rhoGn))} \to \juerg{(\mui,(\rhoi,\rhoGi))}
  \quad \hbox{\pier{strongly} in $\C0\Vp\times\C0\calH$}\,,
  \non
\Eeq
and we infer that
\begin{align*}
& \juerg{ (\mu,(\rho,\rhoG))}(\tn) = \juerg{(\mun,(\rhon,\rhoGn))}(0) \to 
\juerg{(\mui,(\rhoi,\rhoGi))}(0) = \juerg{(\mus,(\rhos,\rhoGs))}\\[1mm]
 &  \hbox{weakly in $\Vp\times\calH$}.
\end{align*}
By comparing with~\eqref{toomega}, we conclude that $(\mus,\rhos,\rhoGs)=(\muo,\rhoo,\rhoGo)$,
and the proof is complete.\QED


\section*{Acknowledgments}
\pier{PC and GG gratefully acknowledge some financial support from 
the MIUR-PRIN Grant 2015PA5MP7 ``Calculus of Variations'', 
the GNAMPA (Gruppo Nazionale per l'Analisi Matematica, 
la Probabilit\`a e le loro Applicazioni) of INdAM (Isti\-tuto 
Nazionale di Alta Matematica) and the IMATI -- C.N.R. Pavia.}


\vspace{3truemm}

%
%

{\small%

\Begin{thebibliography}{10}


%
\bibitem{CahH} 
J.W. Cahn, J.E. Hilliard, 
{\em Free energy of a nonuniform system I. Interfacial free energy}, 
J. Chem. Phys. {\bf 2} (1958), 258--267.
\bibitem{Calcol} L. Calatroni, P. Colli,
{\it Global solution to the Allen--Cahn equation with singular potentials and
dynamic boundary conditions}, Nonlinear Anal. {\bf 79} (2013), 12--27. 
\bibitem{CGW} 
C. Cavaterra, M. Grasselli, H. Wu, {\em
Non-isothermal viscous Cahn--Hilliard equation with inertial term and
dynamic boundary conditions},
Commun. Pure Appl. Anal.  {\bf 13}  (2014), 1855--1890.
\bibitem{ChGaMi} L.\ Cherfils, S.\ Gatti, A.\ Miranville,
    {\em A variational approach to a {C}ahn-{H}illiard model 
    in a domain with nonpermeable walls},
    \newblock J.\ Math.\ Sci.\ (N.Y.) {\bf 189} (2013), 604--636.
{\bibitem{ChFaPr} 
R. Chill, E. Fa\v sangov\'a, J. Pr\"uss,
{\em Convergence to steady states of solutions of the
Cahn-Hilliard equation with dynamic boundary conditions},
Math. Nachr. {\bf 279} (2006), 1448--1462.}
\bibitem{CF1} 
	P.\ {C}olli, T.\ {F}ukao, 
	{\em The Allen--Cahn equation with dynamic boundary conditions and mass constraints},
Math. Methods Appl. Sci. {\bf 38} (2015), 3950--3967.
\bibitem{CF2} 
	P.\ {C}olli, T.\ {F}ukao, 
	{\em {C}ahn--{H}illiard equation with dynamic boundary conditions and mass constraint on the boundary}, 
	\newblock J.\ Math.\ Anal.\ Appl. {\bf 429} (2015), 1190--1213.
\bibitem{CF3} 
	P.\ {C}olli, T.\ {F}ukao, 
	{\em Equation and dynamic boundary condition of {C}ahn--{H}illiard type with singular potentials}, 
	\newblock Nonlinear Anal. {\bf 127} (2015), 413--433.
\bibitem{CGKPS}
{P. Colli, G. Gilardi, P. Krej\v{c}\'{\i}, P. Podio-Guidugli, J. Sprekels},
{\em Analysis of a time discretization scheme 
for a nonstandard viscous Cahn--Hilliard system},
ESAIM Math. Model. Numer. Anal. {\bf 48} (2014), 1061--1087.
\bibitem{CGKS1}
{P. Colli, G. Gilardi, P. Krej\v{c}\'{\i}, J. Sprekels}, {\em A vanishing 
diffusion limit in a nonstandard system of phase field equations},
Evol. Equ. Control Theory {\bf 3} (2014), 257--275.
\bibitem{CGKS2}
{P. Colli, G. Gilardi, P. Krej\v{c}\'{\i}, J. Sprekels},
{\em A continuous dependence result for a nonstandard system of phase field equations},
Math. Methods Appl. Sci. {\bf 37} (2014), 1318--1324.

\bibitem{CGPS3}
P. Colli, G. Gilardi, P. Podio-Guidugli, J. Sprekels,
{\it Well-posedness and long-time behaviour for a nonstandard viscous Cahn-Hilliard system},
{SIAM J. Appl. Math.} {\bf 71} (2011) 1849--1870. 
 
\bibitem{CGPS7}
P. Colli, G. Gilardi, P. Podio-Guidugli, J. Sprekels,
{\it Global existence for a strongly coupled Cahn-Hilliard system with viscosity},
{Boll. Unione Mat. Ital. (9)} {\bf 5} (2012), 495--513.

\bibitem{CGPSco}
{P.~Colli, G.~Gilardi, P.~Podio-Guidugli, J.~Sprekels}, {\em Distributed 
optimal control of a nonstandard system of phase field equations},
Contin. Mech. Thermodyn. {\bf 24} (2012), 437--459.

\bibitem{CGPS8}
P. Colli, G. Gilardi, P. Podio-Guidugli, J. Sprekels,
{\it Continuous dependence for a nonstandard Cahn-Hilliard system with nonlinear atom mobility},
{Rend. Sem. Mat. Univ. Pol. Torino} {\bf 70} (2012), 27--52. 

\bibitem{CGPS4}
P. Colli, G. Gilardi, P. Podio-Guidugli, J. Sprekels,
{\it An asymptotic analysis for a nonstandard Cahn-Hilliard system with viscosity},
{Discrete Contin. Dyn. Syst. Ser.~S} {\bf 6} (2013), 353--368. 

\bibitem{CGPS6}
P. Colli, G. Gilardi, P. Podio-Guidugli, J. Sprekels,
{\it Global existence and uniqueness for a singular/degenerate Cahn-Hilliard system with viscosity},
{J.~Differential Equations} {\bf 254} (2013), 4217--4244.

\bibitem{CGSco1}
P.~Colli, G.~Gilardi, J.~Sprekels, {\em Analysis and optimal boundary control 
of a nonstandard system of phase field equations},
Milan J. Math. {\bf 80} (2012), 119--149.
\bibitem{CGS-lom}
P. Colli, G. Gilardi, J. Sprekels, {\em Regularity of the solution to a nonstandard
system of phase field equations}, Rend. Cl. Sci. Mat. Nat. {\bf 147}
(2013), 3--19.
\bibitem{CGS0}
P. Colli, G. Gilardi, J. Sprekels, {\em 
On the Cahn--Hilliard equation with dynamic boundary conditions and a dominating boundary potential}, J. Math. Anal. Appl. {\bf 419} (2014), 972--994.
\bibitem{CGS1}
P. Colli, G. Gilardi, J. Sprekels, {\em 
A boundary control problem for the pure Cahn--Hilliard equation with 
dynamic boundary conditions}, Adv. Nonlinear Anal. {\bf 4} (2015), 311--325. 
\bibitem{CGS2}
P. Colli, G. Gilardi, J. Sprekels, {\em 
A boundary control problem for the viscous Cahn--Hilliard 
equation with dynamic boundary conditions}, 
Appl. Math. Optim. {\bf 73} (2016), 195--225.
\bibitem{CGS3}
P. Colli, G. Gilardi, J. Sprekels, {\em On an application of Tikhonov's fixed
point theorem to a nonlocal Cahn-Hilliard type system modeling phase separation},
J. Differential Equations {\bf 260} (2016), 7940--7964.
\bibitem{CGS4}
P. Colli, G. Gilardi, J. Sprekels, {\em 
Distributed optimal control of a nonstandard nonlocal phase field system},
AIMS Mathematics 1 (2016), 225-260.
\bibitem{CGS5}
P. Colli, G. Gilardi, J. Sprekels, {\em 
Distributed optimal control of a nonstandard 
nonlocal phase field system with double obstacle potential},
Evol. Equ. Control Theory 6 (2017), 35-58. 

\bibitem{CGS10}
P. Colli, G. Gilardi, J. Sprekels,
{\it Global existence for a nonstandard viscous
Cahn-Hilliard system with dynamic
boundary condition}, 
\pier{preprint arXiv:1608.00854 [math.AP] (2016), 
1--29 and to appear in {SIAM J. Math. Anal.}}

\bibitem{CGS11}
P. Colli, G. Gilardi, J. Sprekels,
{\it Optimal boundary control of a nonstandard viscous 
CahnÐHilliard system with dynamic boundary condition},
preprint arXiv:1609.07046 [math.AP] (2016), 1--30. 

\bibitem{CS}
P. Colli, J. Sprekels, {\em Optimal control of an Allen--Cahn equation with singular 
potentials and dynamic boundary condition}, SIAM J. Control Optim. {\bf 53} (2015), 213--234.
\bibitem{CGM1} M.\ Conti, S.\ Gatti, A.\ Miranville,
    {\em Attractors for a {C}aginalp model with a logarithmic potential
              and coupled dynamic boundary conditions},
    \newblock Anal. Appl. (Singap.) {\bf 11} (2013), 1350024, 31 pp.
\bibitem{CGM2} M.\ Conti, S.\ Gatti, A.\ Miranville,
    {\em Multi-component Cahn--Hilliard systems with dynamic boundary
 conditions},
 Nonlinear Anal. Real World Appl.  {\bf 25} (2015), 137--166.     
\bibitem{EllSh} 
C.M. Elliott, S. Zheng, 
{\em On the Cahn--Hilliard equation}, 
Arch. Rational Mech. Anal. {\bf 96} (1986), 339--357.

{\bibitem{FiMD1} 
H. P. Fischer, Ph. Maass, W. Dieterich, 
{\em Novel surface modes in spinodal decomposition,} 
Phys. Rev. Letters {\bf 79} (1997), 893--896.}

{\bibitem{FiMD2} 
H. P. Fischer, Ph. Maass, W. Dieterich, 
{\em Diverging time and length scales of spinodal decomposition modes in thin flows},
Europhys. Letters {\bf 42} (1998), 49--54.}

\bibitem{FG}
E. Fried, M.E. Gurtin, {\em Continuum theory of thermally induced phase transitions
based on an order parameter}, Phys.~D {\bf 68} (1993), 326--343.
\bibitem{GaGu} C.\ G.\ Gal, M.\ Grasselli,
     {\em The non-isothermal Allen-Cahn equation with dynamic boundary conditions}, 
     \newblock Discrete Contin.\ Dyn.\ Syst. {\bf 22} (2008), 1009--1040.
\bibitem{GaWa} C.G.\ Gal, M.\ Warma,		
    {\em Well posedness and the global attractor of some quasi-linear
              parabolic equations with nonlinear dynamic boundary
              conditions},
   \newblock  Differential Integral Equations {\bf 23} (2010), 327--358.
\bibitem{GiMiSc} G.\ Gilardi, A.\ Miranville, G.\ Schimperna, 
     {\em On the Cahn-Hilliard equation with irregular potentials and dynamic boundary 
     conditions}, \newblock Commun.\ Pure.\ Appl.\ Anal. {\bf 8} (2009), 881--912.
\bibitem{GMS10}G.\ {G}ilardi, A.\ {M}iranville, G.\ {S}chimperna, 
	{\em Long-time behavior of the {C}ahn--{H}illiard equation with irregular potentials and dynamic boundary conditions},
	\newblock Chin.\ Ann.\ Math.\ Ser.\ B {\bf 31} (2010), 679--712.    
\bibitem{GoMi} G.R.\ Goldstein, A.\ Miranville,
 {\em A {C}ahn-{H}illiard-{G}urtin model with dynamic boundary
              conditions},
  \newblock Discrete Contin. Dyn. Syst. Ser. S {\bf 6} (2013), 387--400.
\bibitem{GMS11} 
	 G.R.\ {G}oldstein, A.\ {M}iranville, G.\ {S}chimperna, 
	{\em A {C}ahn--{H}illiard model in a domain with non-permeable walls},
	\newblock Phys. D {\bf 240} (2011), 754--766.
\bibitem{G}
M.E. Gurtin, {\em Generalized Ginzburg--Landau and Cahn--Hilliard equations based on a
microforce balance},  Phys.~D {\bf 92} (1996), 178--192.
\bibitem{Heida}
M. Heida, {\em Existence of solutions for two types of generalized versions of the
Cahn--Hilliard equation}, Appl. Math. {\bf 60} (2015), 51--90.
\bibitem{Is} H.\ Israel,
 {\em Long time behavior of an {A}llen-{C}ahn type equation with a
              singular potential and dynamic boundary conditions},
 \newblock J. Appl. Anal. Comput. {\bf 2} (2012), 29--56.%
\bibitem{Li} M. Liero, 
 {\em Passing from bulk to bulk-surface evolution in the
              {A}llen-{C}ahn equation},
   \newblock NoDEA Nonlinear Differential Equations Appl.
  {\bf 20} (2013), 919--942.
\bibitem{MRSS} A.\ Miranville, E.\ Rocca, G.\ Schimperna, A.\ Segatti,
   {\em The {P}enrose-{F}ife phase-field model with coupled dynamic
              boundary conditions},
    \newblock Discrete Contin. Dyn. Syst. {\bf 34} (2014), 4259--4290.

{\bibitem{MiZe}
A. Miranville, S. Zelik,
{\em Exponential attractors for the Cahn-Hilliard equation
with dynamic boundary conditions},
Math. Methods Appl. Sci. {\bf 28} (2005), 709--735.}
   
\bibitem{NC}
A.~Novick-Cohen, {\em On the viscous {C}ahn-{H}illiard equation, in
``Material instabilities in continuum mechanics ({E}dinburgh, 1985--1986)''},
Oxford Sci. Publ., Oxford Univ. Press, New York, 1988, pp.~329--342.
\bibitem{PG} {P.~Podio-Guidugli}, {\em Models of phase segregation 
and diffusion of atomic species on a lattice}, Ric. Mat. 
{\bf 55} (2006), 105--118.

{\bibitem{PrRaZh}
J. Pr\"uss, R. Racke, S. Zheng, 
{\em Maximal regularity
 and asymptotic behavior of solutions for the Cahn-Hilliard equation
 with dynamic boundary conditions},  
Ann. Mat. Pura Appl. (4) {\bf 185} (2006), 627--648.}

{\bibitem{RaZh}
R. Racke, S. Zheng, 
{\em The Cahn-Hilliard equation with dynamic boundary conditions}, 
Adv. Differential Equations {\bf 8} (2003), 83--110.}

\bibitem{Simon}
J. Simon,
{\it Compact sets in the space $L^p(0,T; B)$},
{Ann. Mat. Pura Appl.~(4)\/} 
{\bf 146} (1987), 65--96.

\End{thebibliography}

}

\End{document}


\bibitem{Barbu}
V. Barbu,
``Nonlinear Differential Equations of Monotone Type in Banach Spaces'',
Springer,
London, New York, 2010.

\bibitem{Brezis}
H. Brezis,
``Op\'erateurs maximaux monotones et semi-groupes de contractions
dans les espaces de Hilbert'',
North-Holland Math. Stud.
{\bf 5},
North-Holland,
Amsterdam,
1973. 

\bibitem{DiB}
E. DiBenedetto, ``Degenerate parabolic equations'',
Universitext. Springer-Verlag, New York, 1993.

\bibitem{Lions}
J.-L.~Lions,
``Quelques m\'ethodes de r\'esolution des probl\`emes
aux limites non lin\'eaires'',
Dunod; Gauthier-Villars, Paris, 1969.